\documentclass[12pt,onecolumn]{IEEEtran}

\usepackage[utf8]{inputenc}

\usepackage{amsmath}
\usepackage{amssymb}
\usepackage{mathtools}
\usepackage{xcolor}
\usepackage{enumerate}
\usepackage{tikz}
\usetikzlibrary{arrows, arrows.meta, calc, positioning, quotes, shapes}
\usepackage{color}

\usepackage{amsthm}
\newtheorem{theorem}{Theorem}[section]
\newtheorem{definition}[theorem]{Definition}%
\newtheorem{prop}[theorem]{Proposition}%
\newtheorem{lemma}[theorem]{Lemma}% 
\newtheorem{corollary}[theorem]{Corollary}%

\newtheorem{remark}{Remark}[section]

\newtheorem{condition}{Condition}[section]

\IEEEoverridecommandlockouts

\allowdisplaybreaks

\newcommand{\eps}{\varepsilon}
\newcommand{\Rn}{\mathbb{R}^n}

\newcommand{\cminusr}{\mathbb{C}\setminus\mathbb{R}}

\newcommand{\floor}[1]{\left\lfloor #1 \right\rfloor}
\newcommand{\abs}[1]{\left| #1 \right|}
\newcommand{\parens}[1]{\left( #1 \right)}
\newcommand{\cbrackets}[1]{\left\{#1\right\}}
\newcommand{\sbrackets}[1]{\left[#1\right]}
\newcommand{\foralltext}{\textrm{for all }}

\def\one{\mbox{1\hspace{-5.5pt}\fontsize{9.8}{11}\selectfont\textrm{1}}}
\newcommand{\indicator}[1]{\one_{\cbrackets{#1}}}

\newcommand{\landau}[2]{\mathcal{O}_{#1}\parens{#2}}

\DeclareMathOperator{\tr}{tr}
\newcommand{\trp}[1]{\tr\parens{#1}}
\newcommand{\pnorm}[2]{\left\Vert#1\right\Vert_{#2}}
\newcommand{\inorm}[1]{\pnorm{#1}{\infty}}

\newcommand{\expec}[2][]{E_{#1}\left[#2\right]}
\newcommand{\prob}[2][]{P_{#1}\left(#2\right)}

\newcommand{\qpartial}[2]{u_{#1}^{#2}}
\newcommand{\zoomquantizer}{Q_{K}^{\Delta_t}}
\newcommand{\qfixed}{U_N}
\newcommand{\DeltaN}{\Delta_{(N)}}

\newcommand{\piN}{\pi_N}
\newcommand{\xpi}{x_{*,N}}
\newcommand{\xpiparens}{(\xpi)}
\newcommand{\Deltapi}{\Delta_{*,N}}
\newcommand{\epi}{e_{*,N}}

\newcommand{\abslam}{\abs{\lambda}}
\newcommand{\systemzoom}{(\sqrt{2}\abslam+1)}
\newcommand{\systemzoomnp}{\sqrt{2}\abslam+1}

\newcommand{\Two}{T_{w_0}}
\newcommand{\Tcal}{\mathcal{T}}
\newcommand{\Tcali}[1]{{\Tcal_{#1}}}

\newcommand{\phichain}{\cbrackets{\phi_t}_{t=0}^\infty}
\newcommand{\phispace}{\mathbb{X}}
\newcommand{\phiborel}{\mathcal{B}(\mathbb{X})}

\newcommand{\ourchain}{\cbrackets{(x_t,\Delta_t)}_{t=0}^\infty}

\newcommand{\state}[1]{x_{#1}, \Delta_{#1}}
\newcommand{\statep}[1]{\parens{\state{#1}}}

\newcommand{\SchemeP}{Scheme P$(\beta,\eps)$}

\title{An Asymptotically Optimal Two-Part Fixed-Rate Coding Scheme for
  Networked Control with Unbounded Noise   \thanks{The authors are with the Department of Mathematics and
                Statistics, Queen's University, Kingston, Ontario, Canada, K7L
                3N6.  Email: jonathan.keeler@queensu.ca, tamas.linder@queensu.ca,
                yuksel@queensu.ca. This work was presented in part at
                the 2022 IEEE International Symposium on Information Theory, ISIT 
                2022 \cite{KeLiYu22}. }
}
\author{Jonathan Keeler, Tamás Linder, Serdar Yüksel}

\begin{document}

\maketitle

\begin{abstract}
  It is known that under fixed-rate information constraints, adaptive
  quantizers can be used to stabilize an open-loop-unstable linear
  system on $\mathbb{R}^n$ driven by unbounded noise. These adaptive
  schemes can be designed so that they have near-optimal rate, and the
  resulting system will be stable in the sense of having an invariant
  probability measure, or ergodicity, as well as boundedness of the
  state second moment. Although structural results and information
  theoretic bounds of encoders have  been
  studied, the performance of such adaptive fixed-rate
  quantizers beyond stabilization has not been addressed. In this
  paper, we propose a two-part adaptive (fixed-rate) coding scheme
  that achieves state second moment convergence to the classical
  optimum (i.e., for the fully observed setting) under  mild moment
  conditions on the noise process. The first part, as in prior work, leads to ergodicity (via positive Harris recurrence)
  and the second part ensures that the state second moment converges
  to the classical optimum at high rates. These results are
  established using an intricate analysis which uses random-time
  state-dependent Lyapunov stochastic drift criteria as a core tool.
\end{abstract}


\section{Introduction}\label{sec:intro}
\emph{Networked control} or \emph{information-constrained control}
refers to control systems in which the controllers, sensors, and
systems (actuators/plants) are connected through communication
channels or a data-rate constrained network. Thus, there may be a data
link between the sensors (which collect information), the controllers
(which make decisions), and the actuators (which execute the
controller commands). Moreover, the sensors, controllers and the plant
themselves could be geographically separated. For such
information-constrained control (or networked control) systems, one
needs to jointly design encoders and controllers for satisfactory
performance, which may have stability or optimality as a design
objective.

In cases where stability is the primary design objective, one is
typically concerned with the minimum capacity above which
stabilization is possible, and there are many results of this flavour
in the existing literature (these are discussed in detail  in
Section \ref{sec:lit_review}).

This paper is primarily concerned with optimality as the primary
design objective. In this context, the notion of optimality is the
minimization of some cost function over a specified time horizon
(which here is infinite), and for the kinds of systems we consider
here, one seeks asymptotic bounds on the cost function as the data
rate becomes large.

\subsection{Problem Statement}\label{sec:problem_statement}
Let us first introduce the system to be controlled under no information constraints. The optimal cost in this ``classical'' case will yield a lower bound over all information-constrained policies. Consider the linear (but not necessarily Gaussian) discrete-time multi-dimensional control system,
\begin{equation}\label{eq:openloop_system}
    x_{t+1} = Ax_t + Bu_t + w_t,
\end{equation}
where the state process $x_t$, the control $u_t$, and the noise process $w_t$ live in $\Rn$, $A$ and $B$ are $n\times n$ real matrices, and we assume $B$ to be invertible (this can be relaxed to controllability by a sampling argument for stability \cite{AndrewJohnstonReport}, though in this case the optimality results we present are not maintained).

The initial state $x_0$ may be distributed according to some
probability measure $\nu$ on $\mathbb{R}^n$, as long as $x_0\sim\nu$
admits at least the same finite moments as the noise process
$w_t$. That is,  whenever
$E\big[ {\inorm{w_t}^\beta} \big]  < \infty$ for some $\beta>0$,  we
have $E\big[{\inorm{x_0}^\beta}\big] < \infty$ (this is trivially
satisfied if $\nu = \delta_x$ for some $x\in\mathbb{R}^n$). In
addition, we suppose that the initial state $x_0$ is independent of
the noise process $w_t$.

In the classical, fully observed setting, at each time stage $t\geq 0$, the controller has access to the history $I_t = x_{\sbrackets{0, t}} = \parens{x_0,\ldots,x_t}$. An admissible control policy $\gamma$ is a sequence of Borel measurable mappings $\cbrackets{\gamma_t; t\geq 0}$ where $\gamma_t : I_t \to \Rn$ is such that it produces the control $u_t = \gamma_t(I_t)$ at each time stage.

The noise process $\cbrackets{w_t}_{t=0}^\infty$ is assumed to be
i.i.d., zero-mean with covariance matrix
$\Sigma = \expec[]{w_0w_0^\top}$. Furthermore, we suppose that the noise
process admits a pdf $\eta$ with respect to the Lebesgue measure on
$\Rn$ which is positive everywhere (and thus has unbounded support).

For an $n\times n$ positive definite matrix $Q$, the optimal control problem is to choose a policy $\gamma$ which minimizes the infinite-horizon average quadratic form,
\begin{equation}\label{eq:ergodic_cost}
    \limsup_{T\to\infty} \frac{1}{T}E^\gamma\sbrackets{\sum_{t=0}^{T-1}x_t^{\top}Qx_t},
\end{equation}
where the expectation above is with respect to the policy $\gamma$ and
$x^\top$ denotes the transpose of the column vector $x\in
\mathbb{R}^n$.

The following statement is a special case of a
well-known result
for a more general setup (where, e.g.,  $B$ may not be invertible) and can
be obtained by solving the algebraic Riccati equation (see, e.g.,
\cite[Chap.~4.1]{ber05}. However, we provide a short and direct
proof in the Appendix. 

\begin{prop}\label{prop:fully_observed_optimal}
For the fully observed setup described above, the optimal control policy is $u_t = -B^{-1}Ax_t$, achieving an optimal cost of $\trp{Q\Sigma}$ where $\trp{\cdot}$ is the trace operator.
\end{prop}

In contrast to this idealized fully observed setup, in this paper we
assume that the controller only has access to $x_t$ through a discrete
noiseless channel of capacity $C$ bits. We assume that the encoder is
causal. In particular, letting $\mathcal{M}$ be a finite alphabet of
cardinality $\abs{\mathcal{M}}= \lfloor 2^C\rfloor$, the encoder is specified by a
quantization policy $\Pi$, which, with $\mathbb{X} = \mathbb{R}^n$ as
the state space, is a sequence of functions
$\cbrackets{\eta_t}_{t=0}^\infty$ of the type
$\eta_t : \mathcal{M}^t\times\mathbb{X}^{t+1}\to\mathcal{M}$. At time
$t$, the encoder transmits the $\mathcal{M}$-valued message
\begin{equation*}
    q_t = \eta_t(I^e_t),
\end{equation*}
where $I^e_0 = x_0$ and $I^e_t = (q_{[0,t-1]}, x_{[0,t]})$ for $t \geq 1$. The collection of all such zero-delay policies is called the set of admissible quantization policies and is denoted by ${\bf \Pi}_A$. 

Upon receiving $q_t$, the receiver generates the control $u_t$, also
without delay. A zero-delay controller policy is a sequence of
functions $\gamma = \cbrackets{\gamma_t}_{t=0}^\infty$ such that 
$\gamma_t : \mathcal{M}^{t+1}\to\mathbb{U}$, where $\mathbb{U}=\mathbb{R}^n$ is the
control action space, so that $u_t = \gamma_t(I^d_t)$, where
$I^d_t=q_{[0,t]}$.

\begin{figure}[h]
\centering
\scalebox{1}{
\begin{tikzpicture}[
    node distance = 10mm and 15mm,
    base/.style={rectangle, draw=black, fill=black!5, thick},
    smallblock/.style={base, minimum size=6mm},
    mediumblock/.style={base, minimum width=15mm, minimum height=8mm},
    bigblock/.style={base, minimum width=20mm, minimum height=10mm}]
%Nodes
\node[bigblock, align=center]                             (channel)           {Noiseless\\ Channel};
\node[smallblock, left=of channel]          (coder)             {Coder};
\node[smallblock, right=of channel]         (controller)        {Controller};
\node[mediumblock, below=of channel]        (plant)             {Plant};

%Lines
\draw [-{Latex[length=2mm]}] 
    (coder)         edge ["{$q$}"]      (channel)
    (channel)       edge ["{$q'=q$}"]     (controller);

\draw [-{Latex[length=2mm]}] 
    (controller.east) -- +(0.75,0) |- node[above, pos=0.9] {$u$} (plant);
\draw [-{Latex[length=2mm]}]
    (plant.west) -- node[above, pos=0.2] {$x$} +(-4,0) |- (coder);

\end{tikzpicture}
}
\caption{Block diagram of the communication and control loop.}
\label{fig:blockdiagram}

\end{figure}
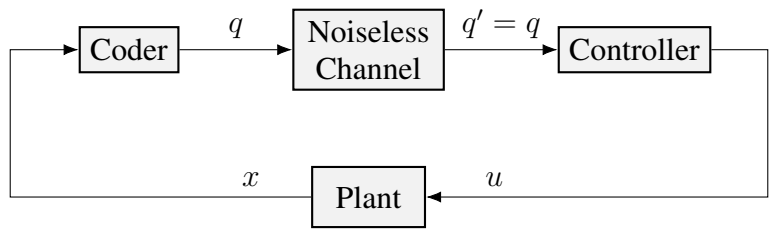

Thus, the data rate is fixed, and we assume zero coding delay. In this setup, it becomes necessary to describe not just a control policy, but also a coding scheme with which to communicate information about the current state vector.

\subsection{Literature Review}\label{sec:lit_review}

For systems of this nature, various authors have obtained the minimum 
channel capacity above which stabilization is possible, under various
assumptions on the system and the admissible coders and
controllers. Here, ``stabilization'' can be in several senses, for
example positive Harris recurrence, asymptotic mean stationarity or
more generally, limiting moment stability of the state.

Such a result is usually referred to as a \emph{data-rate theorem} and
takes the following form, with $\cbrackets{\lambda_i}$ being the
eigenvalues of the system matrix $A$:
\begin{equation}
    C > R_{\textrm{min}} \coloneqq \sum_{\abs{\lambda_i} \geq 1} \log_2\abs{\lambda_i}
\end{equation}
That is, the capacity must exceed the sum of the unstable eigenvalue logarithms.

Some of the earliest works in this context are \cite{Brockett} and \cite{Baillieuil99}. More general versions of the data-rate theorem have been proven in \cite{Tatikonda} and \cite{hespanha2002towards}. For noisy systems and mean-square stabilization, or more generally, moment-stabilization, analogous data-rate theorems have been proven in \cite{NairEvans} and \cite{SahaiParts}, see also \cite{Martins,SavkinSIAM09}.

In \cite{YukTAC2010,YukMeynTAC2010}, a joint fixed-rate coding and control scheme is given which, in the scalar case $n=1$ with unstable eigenvalue $\abslam \geq 1$ and where $w_t$ is Gaussian, stabilizes the system \eqref{eq:openloop_system} while being nearly rate-optimal, in that the rate used satisfies only $C > \log(\abslam+1)$. This is achieved using an adaptive uniform quantization scheme, where the quantizer bin sizes "zoom" in and out exponentially to track the state $x_t$. Here, the notion of stability is ergodicity and finiteness of all limiting system moments. By increasing a sampling period $T$, the achievable rate $\frac{1}{T}\log(\abslam^T+1)$ gets arbitrarily close to $C > R_{\textrm{min}} = \abslam$ \cite[Theorem 2.3]{AndrewJohnstonReport}. This scheme can be generalized to one which stabilizes the multi-dimensional system \eqref{eq:openloop_system} (where the noise is more general than Gaussian) using a similar approach \cite{AndrewJohnstonReport}. Furthermore, this leads to a closed loop system which is positive Harris recurrent (and hence, ergodic) and admits finite limiting system second moment \cite[Theorem 2.2]{AndrewJohnstonReport}. For related recent fixed-rate constructions which also utilize modest delay, we refer the reader to \cite{kostina2021exact} and \cite{sabag2020stabilizing}.

Despite being near rate-optimal for achieving stability (i.e., finite
asymptotic system moments), the schemes in \cite{YukTAC2010,
  YukMeynTAC2010, AndrewJohnstonReport} have not been shown to yield
second moment convergence to the classical optimum $\trp{Q\Sigma}$ as
the data rate $C$ grows large.

In the literature, information theoretic relaxations of the problem
noted have been studied, for obtaining both lower and upper
bounds. Lower bounding methods typically build on replacing the number
of bins with entropy of the quantization symbols, or the latter with
mutual information bounds, and the use of Shannon lower bounding
techniques. Upper bounding methods include entropy coding and
dithering methods; dithering \cite{zamir1992universal} ``uniformizes''
the noise even for low rates (though which critically requires the
presence of common randomness at the encoder and the decoder); see
\cite{Silva1}, \cite{ostergaard2021stabilizing}
\cite{silva2015characterization}, \cite{stavrou2018zero},
\cite{cuvelier2022time}. Using ergodicity and invariance
properties, \cite{cuvelier2022time} has established time-invariant
(though still variable-rate) coding schemes using dithering. Making
use of lattice quantization performance bounds developed in
\cite{kostina2017data},  without the use of of dithering  \cite{kostina2019rate} 
investigated distributed control with lattice quantization followed by entropy coding (and thus
also with variable rate coding) and established near optimality of
such  schemes for high-rates.

In this paper, we will only use fixed-rate codes, which may be more
suitable for a large class of zero-delay systems.

% At high rates one can obtain a uniform error approximation, without
% the presence of randomized dithering. Uniform scalar quantizers, in
% addition to their simplicity, are near-optimal in the limit of
% high-rate distortion \cite{GishPierce}. In the limit of low
% distortion (high rate), \cite{LinderZamir} studied stationary
% encoding of a stable stationary process and showed that memoryless
% quantizer followed by a conditional entropy coder is near
% optimal. This study, however, has assumed stationarity of the source
% where adaptation is not needed. We also refer the reader also to
% \cite{kawanyukselTAC2020} for high-rate results as well as relations
% with metric and topological entropy of dynamical systems.

A further common approach in the literature has been to minimize the
directed information \cite{massey90}  subject to the distortion rate: 
\[
  \limsup_{N \to \infty} \min\Big\{ \frac{1}{N} I(X^N \to \hat{X}^N):
  \frac{1}{N}E\Big[\sum_{k=0}^{N-1} \big(X_k - \hat{X}_k\big)^2\Big]
  \leq D \Big\}.
\]
This leads to an information theoretic lower bound to the optimal
estimation error subject to an information rate constraint. There has
been a surge of research activity on this problem since
\cite{gorbunov1973nonanticipatory}, where explicit solutions, bounds,
as well as convex analytic numerical solutions (including via
semi-definite programming) have been presented; see, e.g., 
\cite{banbas89,TatikondaSahaiMitter,charalambous2014nonanticipative,tanaka2015semidefinite,tanaka2016semidefinite,tanaka2018lqg,stavrou2018zero,derpich2012improved,stavrou2018optimal,stavrou2021asymptotic}. For noisy channels, using channel-source coding separation based methods via the rate-distortion function and Shannon capacity duality can be utilized to establish tightness results, especially for the Gaussian case. When an additive Gaussian channel is present, sequential rate-distortion theoretic ideas presented above in (also via generalizing the scalar Gaussian analysis \cite{banbas89}) lead to explicit optimality conditions studied in \cite{tanaka2016semidefinite} and \cite{stavroutanaka2019time} (see also \cite{kostina2019rate})

%
%A particularly important approach is to minimize the directed information. This leads to an information theoretic lower bound to the optimal estimation error subject to an information rate constraint. There has been a surge of research activity on this problem since \cite{gorbunov1973nonanticipatory} where explicit solutions, bounds, as well as convex analytic numerical solutions (including via semi-definite programming) have been presented in \cite{banbas89, TatikondaSahaiMitter, charalambous2014nonanticipative, tanaka2015semidefinite, tanaka2016semidefinite, tanaka2018lqg, stavrou2018zero, derpich2012improved, stavrou2018optimal, stavrou2021asymptotic} and \cite{stavroutanaka2019time} (see also \cite{kostina2019rate}).
%

%For noisy channels, channel-source coding separation based
%methods via the rate-distortion function and Shannon capacity duality
%can be utilized to establish tightness results, especially for the
%Gaussian case.  When an
%additive Gaussian channel is present, sequential rate-distortion
%theoretic ideas  (also via generalizing the scalar
%Gaussian analysis \cite{banbas89}) lead to explicit optimality
%conditions  in \cite{tanaka2016semidefinite} and
%\cite{stavroutanaka2019time} (see also
%\cite{kostina2019rate}). %These may find important operational applicability when sufficient conditions on the system and channel models (at least numerically) on the optimality of linear policies can be obtained.

We refer to \cite{YukselLQGQuantizationTAC} and \cite{yuksel2019note}
for a detailed review on structural results for optimal coding of
controlled linear systems. Infinite horizon zero-delay coding for
linear systems is studied in \cite{YukLinZeroDelay, ghomi2021zero}.

\subsection{Contributions}
We study the class of stochastic linear systems having transition
dynamics given in \eqref{eq:openloop_system} and  driven by unbounded
noise that are to be controlled across a
discrete noiseless channel of finite capacity. For such systems, we
present a novel two-stage coding scheme, in which the first stage is
time-adaptive and stabilizing (in the sense of positive Harris
recurrence and finite limiting system moments), while the second stage
is fixed in time. The first coding stage is a variation on the schemes
in \cite{YukTAC2010,YukMeynTAC2010, AndrewJohnstonReport}.

Crucially utilizing the ergodicity results of the first coding stage,
we show that this two-part coding scheme attains convergence of the
limiting system second moment to the classical optimum as the data
rate $C$ grows large, with explicit rate of convergence. While
multi-stage quantization schemes have been studied before in the
source coding literature \cite{gray1982multistage_vq}, our
implementation is novel in that one stage of the code is time-adaptive
(though still with fixed-rate) and stabilizing, and the second stage
is designed for near-optimal performance at high-rates. An inspiration
for this approach also comes from Berger \cite{BergerWiener} and Sahai
\cite{SahaiISIT04}. Our analysis complements recent studies in the
literature \cite{kostina2019rate} \cite{cuvelier2022time} which have
utilized variable rate codes (without and with dithers, respectively),
but more importantly, it presents a two-part architecture in
code design, separately targeting stability and optimality.

To our knowledge, this is the first scheme proven to have this
convergence property for systems of this type (in particular, with
unbounded noise) in networked control. The algorithm to be presented
is relatively simple, and being fixed-rate, highly practical, but its
performance analysis requires quite involved technical steps which
necessitate a technically involved  analysis using, among other tools, random-time
state-dependent stochastic Lyapunov drift conditions \cite{YukTAC2010}
and a careful analysis of the performance of finite-level uniform
quantizers for unbounded sources that satisfy certain moment
constraints. 

\subsection{Organization}
The paper is organized as follows. Section \ref{sec:prelims} contains relevant preliminaries, some background on stochastic stability for general state-space Markov chains, and statements of some useful random-time Lyapunov drift theorems that are central to our analysis.

Section \ref{sec:scalar} presents the construction of our two-part
scheme in the scalar case for simplicity, and provides a somewhat
detailed proof program to guide the reader. Section \ref{sec:vector}
provides the construction of our two-part scheme in the vector case,
in full generality. A detailed proof program in the vector case is
provided in Section \ref{sec:proof_program_vector}.

%Section \ref{sec:simulation} presents a simulation that illustrates our two-part scheme and main results.

Finally, as many of the proofs are rather mechanical and quite
involved, the Appendix contains complete proofs of many key results
stated in the main body of the paper.

\section{Preliminaries}\label{sec:prelims}

\subsection{Definitions and Conventions}\label{ssec:definitions}
We denote the nonnegative and strictly positive reals by $\mathbb{R}_+$ and $\mathbb{R}_{++}$, respectively. We let $\mathbb{N}$ denote the nonnegative integers.

For $x\in\Rn$ and $p\in[1,\infty)$ we denote $\pnorm{x}{p} \coloneqq \parens{\sum_{i=1}^n \abs{x^i}^p }^{\frac{1}{p}}$. For $p=\infty$ we denote $\inorm{x} \coloneqq \max_{1\leq i\leq n} |x^i|$.

It is well-known that $\pnorm{\cdot}{p}$ is a norm on $\Rn$ for all $p\in[1,\infty]$. The following pair of inequalities is also well-known (e.g., see \cite[Exercise 3.5(a)]{BigRudin}):
\begin{prop}
Suppose $1\leq p \leq q \leq \infty$. Then for any $x\in\Rn$,
\begin{equation}\label{eq:pnorm_inequality}
    \pnorm{x}{q} \leq \pnorm{x}{p} \leq
    n^{\frac{1}{p}-\frac{1}{q}}\pnorm{x}{q}, 
\end{equation}
where $\frac{1}{\infty}$ is taken to be $0$ by convention.
\end{prop}
The first inequality follows by a simple renormalization argument,
while the second follows from H\"{o}lder's inequality.

For a matrix $V\in\mathbb{R}^{m\times n}$, we will refer to the $(i,j)$-th component either as $V_{ij}$ or as $\sbrackets{V}_{ij}$. The latter notation will be used when $V$ takes on an expression involving square brackets (e.g., an expectation) so as to avoid ambiguity.

For a matrix $A\in\mathbb{R}^{m\times n}$ we define its $\infty$-norm as $\inorm{A} \coloneqq \max_{1\leq i\leq m} \sum_{j=1}^n \abs{A_{ij}}$. This norm is consistent with the vector $\infty$-norm in the following sense. Let $v\in\Rn$ and $A\in\mathbb{R}^{m\times n}$. Then,
\begin{equation}\label{eq:matrixvector_inequality}
    \inorm{Av} \leq \inorm{A}\inorm{v}.
\end{equation}
We will find it useful to use the Landau notation for comparing function asymptotics. For two functions $f,g:[a,\infty)\to\mathbb{R}_+$, where $a\in\mathbb{R}$, we say that $f = \landau{u}{g}$ if for all $u$ sufficiently large one has $f(u) \leq cg(u)$ for some constant $c>0$. This is equivalent to the condition that $\limsup_{u\to\infty}\frac{f(u)}{g(u)} < \infty$.

Finally, for a vector-valued random variable $X$ we denote its ``tail function'' as
\begin{equation*}
    T_X(u) \coloneqq \prob[]{\inorm{X} > u},\quad u\geq 0.
\end{equation*}

\subsection{Stochastic Stability}\label{sec:stability_prelims}
In this section we provide some brief background on stochastic stability, particularly that which will be relevant to the stabilizing properties of the two-stage scheme we present (namely, positive Harris recurrence).

Suppose $\phichain$ is a time-homogeneous Markov chain with state space $\phispace$, where $\phispace$ is a complete separable metric space that is locally compact; its Borel sigma algebra is denoted $\phiborel$. The transition probability is denoted by $P$, so that for any $\phi\in\phispace$ and $A\in\phiborel$, the probability of moving in one step from state $\phi$ to the set $A$ is given by $\prob[]{\phi_{t+1}\in A\ |\ \phi_t = \phi} \eqqcolon P(\phi,A)$. For any $n\geq 2$, the $n$-step transitions $\prob[]{\phi_{t+n} \in A\ |\ \phi_t=\phi} \eqqcolon P^n(\phi,A)$ are obtained recursively in the usual way:
\begin{equation*}
    P^n(\phi, A) = \int_\mathbb{X}P(y,A)P^{n-1}(\phi, dy).
\end{equation*}
The transition law acts on measurable functions $f:\phispace\to\mathbb{R}$ and measures $\mu$ on $\phiborel$ via
\begin{equation*}
    Pf(\phi) \coloneqq \int_\phispace P(\phi, dy)f(y)\ = \expec[]{f(\phi_{t+1})\ |\ \phi_t=\phi} ,\ \ \foralltext\phi\in\phispace
\end{equation*}
and
\begin{equation*}
    \mu P(A) \coloneqq \int_\phispace\mu(d\phi)P(\phi, A),\ \ \foralltext A\in\phiborel.
\end{equation*}
A probability measure $\pi$ on $\phiborel$ is called invariant if $\pi P = \pi$, that is,
\begin{equation*}
    \int_\phispace \pi(d\phi)P(\phi, A) = \pi(A),\ \ \foralltext A\in\phiborel.
\end{equation*}
For any initial probability measure $\nu$ on $\phiborel$ we can construct a stochastic process with transition law $P$ and $\phi_0 \sim \nu$. We let $P_\nu$ denote the resulting probability measure on the sample space, with the usual convention that $\nu = \delta_\phi$ (i.e., $\nu(\{\phi\})=1$) when the initial state is $\phi \in \phispace$. When $\nu = \pi$ is invariant, the resulting process is stationary.

There is at most one stationary solution under the following irreducibility assumption. For a set $A\in\phiborel$ we denote,
\begin{equation}
    \tau_A \coloneqq \min\cbrackets{t\geq 1:\phi_t\in A}.
\end{equation}

% definition formatting
\begin{definition}
Let $\varphi$ denote a $\sigma$-finite measure on $\phiborel$.
\end{definition}
\begin{enumerate}[\hspace{1em}(i)]
    \item The Markov chain is called \textit{$\varphi$-irreducible} if for any $\phi\in\phispace$ and $B\in\phiborel$ satisfying $\varphi(B) > 0$ we have
    \begin{equation*}
        \prob[\phi]{\tau_B < \infty} > 0.
    \end{equation*}
\end{enumerate}
\begin{enumerate}[\hspace{0em}(ii)]
    \item A $\varphi$-irreducible Markov chain is \emph{aperiodic} if for any $\phi\in\phispace$ and any $B\in\phiborel$ satisfying $\varphi(B) > 0$, there exists $n_0 = n_0(\phi, B)$ such that for all $n\geq n_0$,
    \begin{equation*}
        P^n(\phi, B) > 0.
    \end{equation*}
\end{enumerate}
\begin{enumerate}[\hspace{-1.05em}(iii)]
    \item A $\varphi$-irreducible Markov chain is \textit{Harris recurrent} if $\prob[\phi]{\tau_B < \infty} = 1$ for any $\phi\in\phispace$ and any $B\in\phiborel$ satisfying $\varphi(B) > 0$. It is \textit{positive Harris recurrent} if in addition there is an invariant probability measure $\pi$.
\end{enumerate}\vspace{0.6em}

The notion of full and absorbing sets will be useful to us.

\begin{definition}
For a $\varphi$-irreducible Markov chain $\phichain$, a set $A\in\phiborel$ is called full if $\varphi(A^C) = 0$.
\end{definition}
\begin{definition}
A set $A\in\phiborel$ is called absorbing if $P(x,A) = 1$ for all $x\in A$.
\end{definition}

Finally, we define the notion of small sets for a Markov chain.
\begin{definition}
A set $C\in\phiborel$ is $(m,\delta,\nu)$-small on $(\phispace,\phiborel)$ (for integer $m\geq 1$, $\delta \in (0,1]$ and a probability measure $\nu$ on $(\phispace,\phiborel)$) if for all $x\in C$ and $B\in\phiborel$,
\begin{equation*}
    P^m(x,B) \geq \delta\nu(B).
\end{equation*}
A set is called small (or sometimes $m$-small) if it is $(m,\delta,\nu)$-small for some $(m,\delta,\nu)$.
\end{definition}
Briefly we provide an intuition for small sets. Suppose $C\in\phiborel$ is $(m,\delta,\nu)$-small. Whenever the state $\phi_t$ happens to visit the set $C$, it forgets its entire past with probability at least $\delta>0$ and transitions according to the probability measure $\nu$ over the next $m$ time stages. In this way, the small set $C$ acts as a ``regenerative set'' from which the process can forget its history. This line of investigation leads one to Nummelin's splitting technique \cite{Nummelin1978, AthreyaNey}, which is a key tool in many stability results for irreducible Markov chains.

\subsection{Random-Time State-Dependent Stochastic Lyapunov Drift Conditions for Stability}\label{sec:lyapunov}
In this section, we present a drift condition \cite{YukMeynTAC2010} and use it to show two stability results which will be crucial to our analysis. As in the previous section, we consider a general Markov chain $\phichain$. We consider a sequence of stopping times $\cbrackets{\Tcal_z}_{z=0}^\infty$, measurable with respect to the natural filtration of $\phichain$, such that $\cbrackets{\Tcal_z}_{z=0}^\infty$ is strictly increasing and $\Tcal_0 = 0$. Finally, we will make use of the filtration $\cbrackets{\mathcal{F}_{\Tcal_z}}_{z=0}^\infty$ which is informally the filtration of ``information generated by $\phichain$ up to time $\Tcal_z$'' (for full details,  see \cite{keelerThesis}).

The following is a condition on the general Markov chain $\phichain$ introduced in \cite{YukMeynTAC2010}.

\begin{condition}[Random-Time Lyapunov Drift]\label{cond:lyapunov_random_base}
For a measurable function $V:\phispace\to(0,\infty)$, measurable functions $f,d:\phispace\to[0,\infty)$, a constant $b\in\mathbb{R}$ and a set $C\in\phiborel$, we say that $\phichain$ satisfies the random-time Lyapunov drift condition at $\phi\in\phispace$ if for all $z=0,1,2...$
\begin{equation*}
    \expec[]{V(\phi_\Tcali{z+1})\ |\ \mathcal{F}_\Tcali{z}} \leq V(\phi_\Tcali{z}) - d(\phi_\Tcali{z}) + b\indicator{\phi_\Tcali{z}\in C},
\end{equation*}
and
\begin{equation}\label{eq:driftcriteria}
    \expec[]{\sum_{t=\Tcali{z}}^{\Tcali{z+1}-1}f(\phi_t)\ \big|\ \mathcal{F}_\Tcali{z}} \leq d(\phi_\Tcali{z}).
\end{equation}
when $\phi_0 = \phi$.
\end{condition}

\begin{remark}
Suppose that the stopping times $\Tcali{z}$ are the sequential return times to some set $\Lambda\in\phiborel$, that is $\Tcali{0}=0$ and
\begin{equation*}
    \Tcali{z+1} = \min\cbrackets{t > \Tcali{z} : \phi_t\in\Lambda}.
\end{equation*}
In this case, if one is able to verify for all $\phi\in\Lambda$ that
\begin{equation*}
    \expec[\phi]{V(\phi_{\tau_\Lambda})} \leq V(\phi) - d(\phi) + b\indicator{\phi\in C},
\end{equation*}
and
\begin{equation}\label{eq:driftcriteria_foraset}
    \expec[\phi]{\sum_{t=0}^{\tau_\Lambda-1}f(\phi_t)} \leq d(\phi), 
\end{equation}
then it follows automatically that Condition \ref{cond:lyapunov_random_base} holds at every $\phi\in \Lambda$. Notably, if $C\subseteq\Lambda$ then Condition \ref{cond:lyapunov_random_base} holds at every $\phi\in C$.
\end{remark}

This drift condition, in combination with different assumptions on the functions and sets involved, can lead to many useful results on stability (e.g., \cite[Theorem 2.1]{YukMeynTAC2010}). We present two such results here.

\begin{remark}
The results presented here are variations of \cite[Theorem 2.1]{YukMeynTAC2010}, presented in the form most useful for our application. The proofs of these results draw heavily from the proof program in \cite{YukMeynTAC2010} and rely on supermartingale arguments. For brevity we have omitted the proofs here; for details,  see \cite{keelerThesis} and \cite{YukMeynTAC2010}.
\end{remark}

\begin{lemma}\label{lem:drift_phr}
Suppose $\phichain$ is $\varphi$-irreducible and satisfies Condition \ref{cond:lyapunov_random_base} at all $\phi\in C$, with the restrictions that $C$ is a small set, $\sup_{\phi\in C}V(\phi) < \infty$, $f\equiv 1$ and $d(\phi)\geq 1$. Then the set
\begin{equation*}
    \mathcal{X} \coloneqq \cbrackets{ \phi\in\phispace : \prob[\phi]{\tau_C < \infty} = 1 }
\end{equation*}
is full and absorbing, and the restriction of $\phichain$ to $\mathcal{X}$ is positive Harris recurrent.

If in addition one can show that $\prob[\phi]{\tau_C<\infty} = 1$ for all $\phi\in\phispace$ (i.e., $\mathcal{X} = \phispace$), then $\phichain$ is positive Harris recurrent.
\end{lemma}

\begin{lemma}\label{lem:drift_moment}
Suppose that $\phichain$ is positive Harris recurrent with invariant measure $\pi$ and satisfies Condition \ref{cond:lyapunov_random_base} at some $\phi\in\phispace$. Then,
\begin{equation*}
    \expec[\pi]{f(\phi_t)} \leq b
\end{equation*}
and for any function $g : \phispace \to [0,\infty)$ which is bounded by $f$ in that $g(\cdot)\leq cf(\cdot)$ for some constant $c>0$, we have the following ergodic theorem for $g$,
\begin{equation}\label{eq:driftmomentergodic}
    \lim_{n\to\infty}\frac{1}{n}\expec[\phi_0]{\sum_{t=0}^{n-1}g(\phi_t)} = \expec[\pi]{g(\phi_t)}
\end{equation}
for every $\phi_0\in\phispace$.
\end{lemma}

%%%%%%%%%%%%%%%%%%%%%%%%%%%%%%%%%%%
\section{Scalar Linear Systems}\label{sec:scalar}

We begin with considering the scalar case for simplicity. Here we
consider control of the scalar system,
\begin{equation}
    x_{t+1} = ax_t + bu_t + w_t, 
\end{equation}
where $b\neq 0$ and the noise process $w_t$ is assumed to be i.i.d.,
zero-mean with finite second moment $\sigma^2 =
E\big[{w_0^2}\big]$. Furthermore, we suppose the noise process admits a
pdf $\eta$ with respect to Lebesgue measure on $\mathbb{R}$ which is
positive everywhere. Here we will also suppose that $|a| \geq 1$ so
that the system is open-loop-unstable. That is, if we set  $u_t=0$ for
all $t\ge 1$,
then the system is transient and $x_t$ tends to infinity in magnitude (almost surely and in mean-square).

As in Section \ref{sec:problem_statement}, $x_0$ may be distributed
with some probability measure $\nu$ on $\mathbb{R}$, so long as
$x_0\sim\nu$ admits the same finite absolute moments as the noise process
$w_t$. In addition, we suppose that the
initial state $x_0$ is independent of the noise process $w_t$.

Then the optimal control problem we study is to choose a policy $\gamma$ which minimizes the infinite-horizon quadratic cost,
\begin{equation}\label{eq:scalar_ergodic_cost}
    \limsup_{T\to\infty}\frac{1}{T}E^\gamma\left[{\sum_{t=0}^{T-1}x_t^2}\right].
\end{equation}
By Proposition \ref{prop:fully_observed_optimal}, in the fully
observed setup the optimal control policy is $u_t = -\frac{a}{b}x_t$
which achieves an optimal cost of $\sigma^2$. We suppose now that the
controller only has access to $x_t$ through a discrete noiseless
channel of capacity $C$ bits.

Our goal here is to minimize
$\limsup_{T\to\infty}\frac{1}{T}E^\gamma\big[\sum_{t=0}^{T-1}x_t^2\big]$
which we know is bounded below by $\sigma^2$. As the state will only
be partially observed by the controller, we seek to arrive at
asymptotic bounds on the ``optimality gap''
$\limsup_{T\to\infty}E^\gamma\big[\sum_{t=0}^{T-1}x_t^2\big] -
\sigma^2$ in terms of the channel capacity $C$.

%We now make this precise. As the state is observed through a channel with finite capacity $C$, it is necessary to specify a coding scheme alongside a control policy $\gamma$ (which now must depend on $C$ via the coding scheme). We will do this jointly and in particular, we are interested in the limiting ``high-rate'' case $C\to\infty$.

We say that $\phi_C$ is a joint coding and control scheme for $C>0$ if $\phi_C$ specifies both a coding scheme for communicating over the channel of capacity $C$ as well as a corresponding control scheme at the channel receiver.

The following result demonstrates a limit on achievable rates of convergence.

\begin{lemma}\label{lem:ultimaterate_scalar}
Under any joint coding and control scheme $\phi_C$ we have
\begin{equation}\label{eq:scalar_capacity_bound}
    \limsup_{T\to\infty}\frac{1}{T}\expec[]{\sum_{t=0}^{T-1}x_t^2} - \sigma^2 \geq \frac{a^2\sigma^2}{2^{2C} - a^2}.
\end{equation}
%Therefore, for any ``convergence rate'' function $r:\mathbb{R}_+\to\mathbb{R}_+$ such that the optimality gap is $\landau{C}%{r(C)}$ it must be the case that $\limsup_{C\to\infty}2^{2C}r(C) > 0$.
\end{lemma}
The proof follows the core arguments of the proof of \cite[Theorem 11.3.2]{YukselBasarBook} (with several critical changes, as the assumed limit of distortions utilized in \cite[Theorem 11.3.2]{YukselBasarBook} is to be replaced with a limit of average distortions, requiring additional steps) and is provided in the Appendix.

Intuitively, the above lemma implies that the fastest rate of convergence (of the optimality gap to zero) one can hope for is $2^{-2C}$. Motivated by this, we normalize potential rate functions by $2^{2C}$ and make the following definition.
\begin{definition}\label{def:convergence_scalar}
A joint coding and control scheme $\phi_C$ achieves second moment convergence with rate function $r:\mathbb{R}_+\to\mathbb{R}_+$ if,
\begin{equation*}
    \limsup_{T\to\infty}\frac{1}{T}\expec[]{\sum_{t=0}^{T-1}x_t^2} - \sigma^2 = \landau{C}{\frac{r(C)}{2^{2C}}}.
\end{equation*}
For brevity we may say simply that $\phi_C$ achieves the rate function $r(C)$.
\end{definition}
Intuitively, one seeks to achieve a rate function $r(C)$ which grows slowly in order to maximize the rate of convergence. Lemma \ref{lem:ultimaterate_scalar} implies that the best achievable rate function is the constant function, so any achievable rate function satisfies $\limsup_{C\to\infty}r(C) > 0$.

We will impose the following mild condition on the noise process.
\begin{condition}\label{cond:beta_scalar}
For some $\beta>2$, $\expec[]{|w_0|^\beta} < \infty$. That is, the noise process has finite $\beta$th moment.
\end{condition}

In the scalar case, our main result is the following.
\begin{theorem}\label{thm:ultimate_scalar}
Supposing Condition \ref{cond:beta_scalar} holds with $\beta>2$, for any $\eps \in (0,\beta-2)$ there exists a joint coding and control scheme, denoted \SchemeP, which achieves the exponential rate function $r(C) = 2^{\frac{4}{\beta-\eps}C}$.
\end{theorem}
Intuitively, with only the condition that the noise admits finite
$\beta$th moment, we are able to construct schemes that nearly (as
$\eps\to 0$) achieve convergence of the optimality gap like
$2^{(-2+\frac{4}{\beta})C}$. If $\beta$ is quite large, then the
convergence becomes much closer to $2^{-2C}$ and, in the extreme case
where $w_t$ admits finite moments of all orders (such as in the case
of a Gaussian), one can construct schemes achieving convergence of the
optimality gap at a speed  $\landau{C}{2^{(-2+\delta)C}}$ for any $\delta>0$.

The rest of Section \ref{sec:scalar} is dedicated to the construction
of \SchemeP\ and a high-level proof program of Theorem
\ref{thm:ultimate_scalar}. All proofs are relegated to the Appendix. 

\subsection{Two-Part Code with Uniform Quantization}\label{ssec:scalar_quant}
In this section we describe the joint coding and control scheme used
in proving Theorem~\ref{thm:ultimate_scalar}.

The coding scheme is in two parts where the first part is adaptive in time and the second is fixed. The adaptive part will yield stability, and the fixed part will yield an optimal rate of convergence via simple iterated expectation arguments.

To communicate over a finite capacity channel, will employ uniform
quantizers. Let $M\geq 2$ be an even integer and $\Delta>0$ be a scalar ``bin size''. We define the \emph{scalar modified uniform quantizer} $Q_M^\Delta$ by,
\begin{equation}\label{eq:quantizer_scalar}
    Q_M^\Delta(x) =
    \begin{cases}
        \Delta\floor{\frac{x}{\Delta}} + \frac{\Delta}{2},& \text{if } x\in \left[-\frac{M}{2}\Delta, \frac{M}{2}\Delta \right)\\
        \frac{M}{2}\Delta - \frac{\Delta}{2},& \text{if } x=\frac{M}{2}\Delta\\
        0,& \text{if } |x|>\frac{M}{2}\Delta.
    \end{cases}
\end{equation}
This quantizer uniformly quantizes $x \in \sbrackets{-\frac{M}{2}\Delta, \frac{M}{2}\Delta}$ into $M$ bins of size $\Delta$ and maps all larger $x$ to zero. This requires $M+1$ output levels.

We will use this quantizer for two different purposes.
\begin{itemize}
    \item[(i)] The first is to use adaptive bin sizes which vary with time to achieve stability (in the sense of positive Harris recurrence and finite system moments). Let $K\geq 2$ be an even integer, and suppose that $\cbrackets{\Delta_t}_{t=0}^\infty$ is a sequence of strictly positive ``bin sizes'' varying with time. We will make use of the quantizer $\zoomquantizer$.
    
    \item[(ii)] Secondly, we will use this quantizer to achieve optimal convergence. For a given even number of bins $N\geq 2$, let $\DeltaN$ be a bin size which is a function of $N$. We will make use of the quantizer $U_N^{\DeltaN}$. For brevity, we denote $\qfixed \coloneqq U_N^{\DeltaN}$ and where necessary, specify the dependence of $\DeltaN$ on $N$. Note that this quantizer is fixed in time, in contrast to $\zoomquantizer$.
\end{itemize}

Suppose that the sequence of bin sizes
$\cbrackets{\Delta_t}_{t=0}^\infty$ is such that $\Delta_{t+1}$ is a
function of only $\Delta_t$ and the indicator random variable
$\indicator{|x_t| \leq \frac{K}{2}\Delta_t}$. Also assume that both
the encoder and decoder (controller) know $\Delta_0$. If 
$\zoomquantizer(x_t)$ is sent over the channel, it is possible to
synchronize knowledge of $\Delta_t$ between the quantizer and the
controller since $|x_t|\leq \frac{K}{2}\Delta_t$ if and only if
$\zoomquantizer(x_t)\neq 0$.

\begin{figure}[h]
\centering

\begin{tikzpicture}[
    base/.style={draw=black, fill=black!5, thick},
    arithmetic/.style={base, circle, minimum size=6mm},
    quantizer/.style={base, rectangle, minimum size=6mm},
    channel/.style={base, rectangle, dashed, minimum width=12mm, minimum height=16mm},
    >=Latex]
% Nodes
\node[channel] (Channel) {Channel};
% input of channel nodes
\node[quantizer, left=of Channel, yshift=4mm] (UN) {$U_N^{\DeltaN}$};
\node[arithmetic, left=of UN] (etArith) {\large +};
\node[quantizer, below=of etArith] (QK) {$\zoomquantizer$};
\node[below=of QK] (xt) {$x_t$};
\coordinate[above=of xt, yshift=-5mm] (xtBranch);
\coordinate[above=of QK, yshift=-3.5mm] (QKBranch);
% output of channel nodes
\node[arithmetic, right=of Channel, yshift=-4mm] (xtildeArith) {\large +};
\node[arithmetic, below=of xtildeArith] (facArith) {$\times$};
\node[left=of facArith] (linfac) {$-\frac{a}{b}$};
\node[below=of facArith] (ut) {$u_t$};

% input edges
\draw[->] (xt) -- (QK);
\draw[->] (xtBranch) -- +(-0.8,0) |- node[above, pos=0.7] {$+$} (etArith);
\draw[->] (QK) -- node[left, pos=0.8] {$-$} (etArith);
\draw[->] (etArith) edge ["{$e_t$}"] (UN);

% channel edges
\draw[-] (UN) -- ([yshift=4mm]Channel.west);
\draw[-] (QKBranch) -- ([yshift=-4mm]Channel.west);

\draw[->] ([yshift=-4mm]Channel.east) -- node[above, pos=0.75] {$+$} (xtildeArith);
\draw[->] ([yshift=4mm]Channel.east) -| node[right, pos=0.7] {$+$} (xtildeArith.north);

% output edges
\draw[->] (xtildeArith) edge ["{$\hat{x}_t$}"] (facArith);
\draw[->] (linfac) -- (facArith);
\draw[->] (facArith) -- (ut);

\end{tikzpicture}

\caption{Block diagram of the two-stage coding and control scheme in the scalar case.}
\label{fig:two_stage_scheme_scalar}

\end{figure}
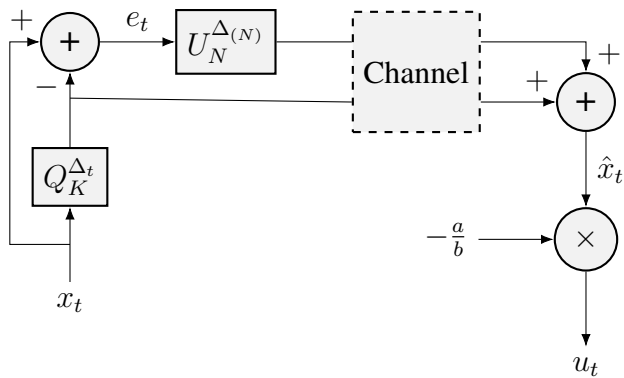

The coding scheme is constructed as follows. For
$\cbrackets{\Delta_t}_{t=0}^\infty$ as above, we calculate the
adaptive quantizer output $\zoomquantizer(x_t)$ and the adaptive
system error $e_t \coloneqq x_t - \zoomquantizer(x_t)$. Then for
integer $N\geq 2$ we use a fixed quantizer $\qfixed$ with bin size
$\DeltaN$ as mentioned above to calculate the fixed quantizer output
$\qfixed(e_t)$. We then send $\zoomquantizer(x_t)$ and $\qfixed(e_t)$
across the noiseless channel where the channel capacity is at least,
\begin{equation}
    C = \log_2{(K+1)} + \log_2{(N+1)}.
\end{equation}
We estimate the state $x_t$ as $\hat{x}_t \coloneqq \zoomquantizer(x_t) + \qfixed(e_t)$. To mirror the fully observed case, the controller applies the control
\begin{equation*}
    u_t = -\frac{a}{b}(\zoomquantizer(x_t) + \qfixed(e_t)) = -\frac{a}{b}\hat{x}_t.
\end{equation*}
The scheme is illustrated in Figure \ref{fig:two_stage_scheme_scalar}.

% Control/coding scheme diagram here

The controlled system dynamics resulting from this scheme are
\begin{align}
    x_{t+1} &= a(x_t - \hat{x}_t) + w_t = a(e_t - \qfixed(e_t)) + w_t.\label{eq:controlled_scalar}
\end{align}

Finally, we describe the adaptive bin size update dynamics where, as in prior work in this context \cite{YukTAC2010,YukMeynTAC2010}, a simple \emph{zooming} scheme is employed. We assume that $K\geq 2$ is even and large enough that $K > |a|$ and choose scalars $\alpha,\rho$ and $L$ such that $\frac{|a|}{K} < \alpha < 1$, $\rho > |a|$ and $L > 0$. We also assume that $\rho \geq K\alpha$. Choose $\Delta_0\geq L$ arbitrarily, then for $t\geq 1$ the bin update is
\begin{equation}\label{eq:binupdate_scalar}
    \Delta_{t+1} = 
    \begin{cases}
        \rho\Delta_t,& \text{if } |x_t| > \tfrac{K}{2}\Delta_t\\
        \alpha\Delta_t,      & \text{if } |x_t|\leq \tfrac{K}{2}\Delta_t,\Delta_t\geq L\\
        \Delta_t,           & \text{if } |x_t|\leq \tfrac{K}{2}\Delta_t,\Delta_t < L.
    \end{cases}
  \end{equation}
The following result is proved, e.g,. in \cite{YukMeynTAC2010}.
\begin{prop}
With dynamics \eqref{eq:controlled_scalar} and \eqref{eq:binupdate_scalar}, the process $\cbrackets{(x_t,\Delta_t)}_{t=0}^\infty$ is a time-homogeneous Markov chain.
\end{prop}
The motivation for this scheme is that the adaptive part leads to stability in the sense of positive Harris recurrence, while the fixed quantizer $\qfixed$ leads to order-optimal convergence of the ergodic second moment \eqref{eq:scalar_ergodic_cost} as the fixed quantization rate $N$ grows large.

The state space for the process $\cbrackets{(x_t,\Delta_t)}_{t=0}^\infty$ highly depends on the following ``countability condition''.
\begin{condition}\label{cond:countability}
There exist relatively prime integers $j,k\geq 1$ such that $\alpha^j\rho^k = 1$. Equivalently, $\log_\alpha\rho$ is rational.
\end{condition}
If this condition holds, then starting from an arbitrary $\Delta_0 > 0$ there exists $\kappa, b\in\mathbb{R}$ such that $\log\Delta_t$ always belongs to a subset of $\mathbb{Z}\kappa + b = \cbrackets{n\kappa+b : n\in\mathbb{Z}}$ (see e.g. \cite[Theorem 3.1]{YukMeynTAC2010}). If the condition fails, then starting from any fixed $\Delta_0$ the set of reachable bin sizes is a dense but countable subset of $\mathbb{R}_{++}$.

We restrict our analysis to the case where Condition \ref{cond:countability} holds. This is not restrictive; it can be shown that for any arbitrary $(\alpha,\rho)$ there exists $(\alpha', \rho')$ arbitrarily close that satisfy Condition \ref{cond:countability}. We let the state space for $\Delta_t$ be
\begin{equation*}
    \Omega_\Delta \coloneqq \cbrackets{\alpha^j\rho^k\Delta_0 : j,k\in\mathbb{Z}_{\geq 0}}.
\end{equation*}
The state space for the Markov chain $\cbrackets{(x_t,\Delta_t)}_{t=0}^\infty$ is then $\mathbb{R}\times\Omega_\Delta$. What remains is to specify additional constraints which complete our proposed scheme.

We assume that Condition \ref{cond:beta_scalar} holds for some $\beta>2$. For any $\eps\in(0,\beta-2)$ we finish our construction of Scheme P$(\beta,\eps)$ by requiring that $\rho > |a|^{\frac{\beta}{\eps}}$ and specifying the dependence of $\DeltaN$ on $N$ as $\DeltaN = 2N^{-1+\frac{2}{\beta-\eps}}$.
\begin{remark}
In this scheme, the constant multiplying $\DeltaN$ is arbitrary for convergence purposes.
\end{remark}
We also impose the following condition in our construction.
\begin{condition}\label{cond:minbinsize_scalar}
The minimum adaptive bin size is at least $\alpha L > \frac{|a|}{K\alpha - |a|}\DeltaN$.
\end{condition}
This may place an implicit dependence of $L$ on the number of fixed quantization bins $N\geq 2$, though if all other parameters remain fixed as $N$ increases then one can ensure Condition \ref{cond:minbinsize_scalar} holds by ensuring that it holds for $N=2$ (since $\DeltaN$ as specified above is monotone decreasing in $N$).

Finally, as $C\to\infty$ we fix $K$ and let $N\to\infty$ to take advantage of fixed quantization results at high rates, to be presented shortly.

Since the proof of Theorem \ref{thm:ultimate_scalar} is rather tedious, we present a somewhat detailed proof program to guide the reader for the scalar setup.

We note that while the proof method for stabilization of our scheme builds on the random-time Lyapunov drift approach introduced in \cite{YukTAC2010,YukMeynTAC2010}, the coupling between the two parts of the coding scheme significantly complicates the analysis. Furthermore, we consider performance bounds as the data rate grows without bound. Altogether, this requires a cautious analysis between moments and high-rate quantization coupled with random-time drift criteria.

\subsection{High-level Proof Program}\label{sec:proof_program_scalar}
In this section we outline the high-level proof program for Theorem \ref{thm:ultimate_scalar}, i.e., the proof that Scheme P$(\beta,\eps)$ achieves the exponential rate function $r(C) = 2^{\frac{4}{\beta-\eps}C}$. Results without complete proofs admit more general sister results in the vector case, which is discussed in detail in Section \ref{sec:proof_program_vector}.

\begin{theorem}
Under \SchemeP\ with $K>|a|$, $\ourchain$ is positive Harris recurrent for every even $N\geq 2$ (i.e., as $C$ grows without bound). Therefore, for every even $N\geq 2$, \SchemeP\ yields a unique invariant measure $\piN$ for the process $\ourchain$.
\end{theorem}
\begin{proof}[Sketch of Proof]
We establish $\varphi$-irreducibility and aperiodicity for $\ourchain$
where $\varphi$ is the product of the Lebesgue and discrete measures
on $\mathbb{R}\times\Omega_\Delta$. The logarithmic function
$V(x,\Delta) = c\log_\alpha\Delta$ is shown to satisfy Condition
\ref{cond:lyapunov_random_base} with $d(x,\Delta)$ constant and
$f\equiv 1$ (i.e., in the form required by Lemma \ref{lem:drift_phr}),
leading to positive Harris recurrence. For a complete proof,  see the
proof of Theorem \ref{thm:PHR} (vector case).
\end{proof}
\begin{remark}
The analysis here is highly similar to that of the proof of \cite[Theorem 3.1]{YukMeynTAC2010}. The only major change is that the upper bound of \cite[Lemma 5.2]{YukMeynTAC2010} for the tail probabilities $\prob[\state{0}]{\tau_\Lambda \geq k}$ must be re-derived in some form, as the out-of-view state dynamics change significantly due to the fixed quantization stage. We address this by providing a similar bound in Lemma \ref{lem:returntime_tailbound} (compare to equation (26) in \cite{YukMeynTAC2010}), which decays suitably fast for summability in the proof program under Condition \ref{cond:beta_scalar}.
\end{remark}

We denote $(\xpi,\Deltapi)\sim\piN$ as the state under invariant measure. This will also induce an invariant measure for the system adaptive error $e_t$, which we denote by $\epi\sim\pi_N^{\textrm{err}}$.

We have the following ergodicity result.
\begin{prop}
The infinite-horizon second moment and the invariant second moment agree, that is
\begin{equation*}
    \lim_{T\to\infty}\frac{1}{T}\expec[]{\sum_{t=0}^{T-1}x_t^2} = \expec[]{(\xpi)^2}.
\end{equation*}
\end{prop}
\begin{proof}[Sketch of Proof]
We are able to show using the drift conditions of Section \ref{sec:lyapunov} that functions $g(x,\Delta)$ that are bounded asymptotically by $\abs{x}^{\beta-\eps}$ satisfy the above ergodicity condition. Since $\eps < \beta - 2$, $g(x,\Delta) = x^2$ is bounded by $\abs{x}^{\beta-\eps}$ and the result follows. For a complete proof, see the proof of Proposition \ref{prop:technical_ergodicity}.
\end{proof}

Finally, we use a simple iterated expectation argument. Suppose that $\statep{0}\sim\piN$. Let $e_0 = x_0 - Q_K^{\Delta_0}(x_0)$ and $x_1 = a(e_0 - \qfixed(e_0)) + Z$, where $Z\sim\eta$ (recall $\eta$ is the distribution of $w_t$). Since we have applied the one-step transition kernel and $\piN$ is the invariant measure, the marginal distributions of $x_0$ and $x_1$ are identical.

For brevity, we denote $s_0 \coloneqq e_0 - \qfixed(e_0)$ so that $x_1 = as_0 + Z$. Supposing that \mbox{$\expec[]{(\xpi)^2} < \infty$} (this will be shown as part of the full proof program), we then have by invariance and iterated expectations that
\begin{align}
    \expec[]{x_0^2} &= \expec[]{x_1^2} = \expec[]{ \expec[]{x_1^2\ |\ s_0} } = \expec[]{ \expec[]{(as_0 + Z)^2\ |\ s_0} }\nonumber\\
    &= \expec[]{a^2s_0^2 + 2as_0\expec[]{Z} + \expec[]{Z^2}} = a^2\expec[]{s_0^2} + \sigma^2\nonumber\\
    &= a^2\expec[]{(e_0 - \qfixed(e_0))^2} + \sigma^2.\label{eqstep:iterated_scalar}
\end{align}
Let $\epi$ denote the system adaptive error under invariant measure, i.e., $\epi = \xpi - Q_K^{\Deltapi}(\xpi)$. Rearranging \eqref{eqstep:iterated_scalar}, we find that the optimality gap is given as
\begin{equation}\label{eq:optimality_gap_scalar}
    \expec[]{\xpiparens^2} - \sigma^2 = a^2\expec[]{(\epi - \qfixed(\epi))^2},
\end{equation}
which is (up to a constant) the distortion of the fixed quantizer $\qfixed$ applied to the random variable $\epi$.

With this in mind, we state the following result for high-rate distortion of $\qfixed$ on sequences of suitably well-behaved random variables. The proof builds on balancing the trade-off between distortion due to the high-rate granular region and the overflow region.
\begin{lemma}\label{lem:MSE_scalar}
Let $\cbrackets{x_N}_{N=2}^\infty$ be a sequence of random variables that satisfy,
\begin{equation}\label{eq:scalar_uniform_bound}
    \sup_{N\geq 2}\expec[]{\abs{x_N}^m} \eqqcolon B_m < \infty
\end{equation}
for some $m > 2$ (not necessarily integer). Set the bin size for the quantizer $\qfixed$ as $\DeltaN = 2N^{-1+\frac{2}{m}}$. Then we have,
\begin{equation*}
    \expec[]{(x_N - \qfixed(x_N))^2} = \landau{N}{N^{-2+\frac{4}{m}}}.
\end{equation*}
\end{lemma}
The proof is mostly mechanical; for a proof, see the proof of
Lemma~\ref{lem:MSE_1} (vector case) in the Appendix.

Briefly, we note that if the sequence $\cbrackets{x_N}_{N=2}^\infty$
is for instance uniformly sub-Gaussian (in the sense that
$\sup_{N\geq 2}E\big[e^{s(x_N)^2}\big] < \infty$ for some $s>0$, then
\eqref{eq:scalar_uniform_bound} holds for every $m>2$ and so it is
possible to achieve distortion asymptotic to
$\mathcal{O}_N(N^{-2+\delta})$ for any $\delta > 0$ by taking $m$
sufficiently large.

However, this sub-Gaussian condition is much stronger than \eqref{eq:scalar_uniform_bound}, and if one follows the proof of Lemma \ref{lem:MSE_scalar} using this stronger tail condition carefully, setting the bin size $\DeltaN$ slightly differently, it is possible to achieve convergence that is faster than $\landau{N}{N^{-2+\delta}}$ (for instance, $\mathcal{O}_N(N^{-2}\ln N)$).

\begin{remark}
  One can  justify that Lemma \ref{lem:MSE_scalar} cannot be
  improved without strengthening the imposed condition or using more
  complex (i.e., non-uniform) quantizers. Consider the special case
  where $x_N = X$ for all $N\geq 2$ and where $X$ admits the
  ``Bucklew-Gallagher'' pdf,
    \begin{equation}\label{eq:bg_pdf}
        p(x) = \frac{1+\delta/2}{(1+|x|)^{3+\delta}}\quad\foralltext x\in\mathbb{R},
    \end{equation}
    for arbitrary $\delta > 0$. This distribution has finite moments
    only of order $m < 2+\delta$ and since $X$ is independent of $N$,
    the problem here is essentially optimal  quantization of 
    $X\sim p$. Let $D_N$ denote the infimum (MSE)
    distortion achievable over all uniform quantizers with $N>0$
    output levels. In \cite{huiNeuhoffUniform}, it is shown  that
    for the source above, $D_N$ satisfies asymptotically (see the
    first equation in \cite[p. 963]{huiNeuhoffUniform}),
    \begin{equation*}
        \lim_{N\to\infty}N^{\frac{2\delta}{2+\delta}}D_N = c_\delta
    \end{equation*}
    for some constant $c_\delta > 0$, so that asymptotically the best distortion achievable by uniform quantizers for the above source $X$ is $D_N = \mathcal{O}_N(N^{-\frac{2\delta}{2+\delta}})$.

    Now we compare this to Lemma \ref{lem:MSE_scalar}. Here, since $X$
    admits finite moments $m$ for any $m < 2+\delta$, by employing our
    uniform quantizer with step size $\DeltaN = 2N^{-1+\frac{2}{m}}$
    we achieve asymptotic distortion
    $\mathcal{O}_N(N^{-2+\frac{4}{m}})$. As we let $m\to 2+\delta$
    from below, $-2+\frac{4}{m}$ becomes arbitrarily close to
    $-\frac{2\delta}{2+\delta}$ and so the result we have stated here
    achieves asymptotic distortion which can be arbitrarily close to
    the best achievable asymptotic distortion one can get using
    uniform quantization. It is in this sense that Lemma
    \ref{lem:MSE_scalar} cannot be improved without strengthening the
    conditions imposed on the sequence $\cbrackets{x_N}$, or without
    employing more complex quantization schemes.
\end{remark}

In light of \eqref{eq:optimality_gap_scalar}, we would like to apply Lemma \ref{lem:MSE_scalar} to the sequence of random variables $\cbrackets{\epi}_{N=2}^\infty$. To do this, we need to establish \eqref{eq:scalar_uniform_bound}. This is ensured by the following result.
\begin{lemma}
Under \SchemeP, the invariant system error has finite $(\beta-\eps)$th moment uniformly in $N\geq 2$. That is,
\begin{equation}\label{eq:sup_error_moments_scalar}
    \sup_{N\geq 2}\expec[]{\abs{\epi}^{\beta-\eps}} < \infty.
\end{equation}
\end{lemma}
\begin{proof}[Sketch of Proof]
With Lyapunov functions $V(x,\Delta)$ and $d(x,\Delta)$ proportional to $\Delta^{\beta-\eps}$ and appropriate set $C$ and constant $b$, we show that these functions satisfy the drift condition \eqref{eq:driftcriteria_foraset}. In particular, we do this for $f$ proportional to $\abs{x}^{\beta-\eps}$ and $f$ proportional to $\Delta^{\beta-\eps}$ which, with the Lyapunov parameters independent of $N$ leads to results of the form \eqref{eq:sup_error_moments_scalar} by Lemma \ref{lem:drift_moment} for the invariant state and adaptive bin size. A simple invariance argument finishes the result. Condition \ref{cond:beta_scalar} is crucial to the proof of this result. For a complete proof, see the proof of Lemma \ref{lem:schemeAmoment} in the Appendix.
\end{proof}
Therefore, in light of \eqref{eq:optimality_gap_scalar} and the above lemma, we find that the optimality gap decays asymptotically at the rate $\mathcal{O}_N(N^{-2+\frac{4}{\beta-\eps}})$. Since $K$ is fixed and $2^C$ is linearly proportional to $N$, we find that the optimality gap is asymptotically in C,
\begin{equation*}
    \expec[]{\xpiparens^2} - \sigma^2 = \landau{C}{\frac{2^{\frac{4}{\beta-\eps}C}}{2^{2C}}},
\end{equation*}
which establishes Theorem \ref{thm:ultimate_scalar}.

\section{The Multi-Dimensional Case}\label{sec:vector}
We now present the more general vector case. Recall that we consider control of the multi-dimensional system \eqref{eq:openloop_system} over a discrete noiseless channel of capacity $C$ bits and the aim is to minimize the infinite-horizon average quadratic form \eqref{eq:ergodic_cost}.

As stated earlier, in the fully observed setup the optimal control policy is $u_t = -B^{-1}Ax_t$ which achieves an optimal cost of $\trp{Q\Sigma}$. In light of this, we seek to arrive at asymptotic bounds on the optimality gap,
\begin{equation*}
\limsup_{T\to\infty}\frac{1}{T}\expec[]{\sum_{t=0}^{T-1}x_t^{\top}Qx_t} - \trp{Q\Sigma}
\end{equation*}
in terms of the channel capacity $C$.

We first generalize Definition \ref{def:convergence_scalar} to the vector case.

\begin{definition}\label{def:convergence_vector}
A joint coding and control scheme $\phi_C$ achieves second moment convergence with rate function $r:\mathbb{R}_+\to\mathbb{R}_+$ if,
\begin{equation*}
    \limsup_{T\to\infty}\frac{1}{T}\expec[]{\sum_{t=0}^{T-1}x_t^{\top}Qx_t} - \trp{Q\Sigma} = \landau{C}{\frac{r(C)}{2^{\frac{2}{n}C}}}.
\end{equation*}
For brevity we may say simply that $\phi_C$ achieves the rate function $r(C)$.
\end{definition}

We refer the reader to \cite[Theorem 4.2]{kawanyukselTAC2020} for a
fundamental lower bound which leads the vector analog to the expression
presented in Lemma \ref{lem:ultimaterate_scalar}. 

\begin{remark}
In the vector case, the  factor of $\frac{1}{n}$ in the above definition comes intuitively from the fact that to tile a hypercube in $\mathbb{R}^n$ of width $L$ it takes $(L/\Delta)^n$ bins of width $\Delta$.
\end{remark}

As in the scalar case, one seeks to achieve a rate function $r(C)$ which grows slowly to maximize the rate of convergence.

We will impose the following condition on the noise process, analogous to Condition \ref{cond:beta_scalar}.
\begin{condition}\label{cond:beta}
For some $\beta > 2$, $\expec[]{\inorm{w_0}^\beta} < \infty$, i.e.,  the noise process has finite $\beta$th moment.
\end{condition}

The following generalization of Theorem \ref{thm:ultimate_scalar} is
our main result.
\begin{theorem}\label{thm:ultimate_result}
Supposing Condition \ref{cond:beta} holds with $\beta > 2$, for any $\eps\in(0,\beta-2)$ there exists a joint coding and control scheme \SchemeP\ which achieves the exponential rate function $r(C) = 2^{\parens{\frac{4}{\beta - \eps}}\frac{1}{n}C}$.
\end{theorem}
The rest of Section \ref{sec:vector} is dedicated to proving Theorem \ref{thm:ultimate_result}.

\subsection{Vector Quantization}\label{sec:vecquant}
In the multi-dimensional case, we will make use of vector quantization. We will use scalar uniform quantizers to define two types of cubic lattice vector quantizers.

Let $M\geq 2$ be an even integer and $\Delta>0$ be a scalar ``bin size''. With $\floor{\cdot}$ the usual floor function, we define the partial uniform quantizer $\qpartial{M}{\Delta} : [-\frac{M}{2}\Delta, \frac{M}{2}\Delta] \to \mathbb{R}$ as
\begin{equation*}
    \qpartial{M}{\Delta}(x) =
    \begin{cases}
        \Delta\floor{\frac{x}{\Delta}} + \frac{\Delta}{2},& \text{if } x\in \left[-\frac{M}{2}\Delta, \frac{M}{2}\Delta \right)\\
        \frac{M}{2}\Delta - \frac{\Delta}{2},& \text{if } x=\frac{M}{2}\Delta.
    \end{cases}
\end{equation*}
Note that $0\not\in \text{range}(u_M^\Delta)$. Then, we define the type I vector quantizer $Q_M^\Delta : \Rn\to\Rn$ as
\begin{equation}
    Q_M^\Delta(x) = 
    \begin{cases}
        \parens{\qpartial{M}{\Delta}(x^i)}_{i=1}^n,& \text{if } \inorm{x} \leq \frac{M}{2}\Delta\\
        0,& \text{if } |x^i| > \frac{M}{2}\Delta^i\ \text{for some } 1\leq i \leq n
    \end{cases}
\end{equation}
and the type II vector quantizer $U_M^\Delta : \Rn\to\Rn$ component-wise as
\begin{equation*}
    \parens{U_M^\Delta(x)}^i = 
    \begin{cases}
        \qpartial{M}{\Delta}(x^i),& \text{if } |x^i| \leq \frac{M}{2}\Delta\\
        0,& \text{otherwise}.
    \end{cases}
\end{equation*}
\begin{remark}
If $n=1$ then the above definitions both correspond to the scalar quantizer \eqref{eq:quantizer_scalar}. For this reason, the scheme we present for the vector case is a direct generalization of the scalar scheme.
\end{remark}

As in the scalar case, we will use these vector quantizers for two different purposes. The type I quantizer will be used with adaptive bin sizes to achieve stability. Let $K\geq 2$ be an even integer, and suppose that $\cbrackets{\Delta_t}_{t=0}^\infty$ is some sequence of strictly positive bin sizes varying with time. We will make use of the quantizer $\zoomquantizer$.

The type II quantizer will be used to achieve convergence to the optimum. For a given
even number of bins $N\geq 2$, let $\DeltaN$ be a bin size which is a
function of $N$. We will make use of the quantizer
$U_N^{\DeltaN}$. For brevity we denote
$\qfixed \coloneqq U_N^{\DeltaN}$ and where necessary, specify the
dependence of $\DeltaN$ on $N$.

\subsection{System in a Jordan Form}\label{sec:sysdesc}

Briefly, we discuss a reduction of the system matrix $A$ into its separate modes. This will have the useful side effect of making $\inorm{A}$ very close to the absolute eigenvalue of each mode.

Let $\cbrackets{\lambda_i}_{i=1}^n$ be the (possibly repeated) eigenvalues of the system matrix $A$. Without loss of generality, we assume that $A$ is in real Jordan normal form, as any matrix $A$ there exists an invertible matrix $P$ such that $\Tilde{A} \coloneqq P^{-1}AP$ is the real Jordan normal form of $A$ \cite[Theorem 3.4.1.5]{MatrixAnalysis}. Let $P\Tilde{x_t} = x_t$ and left-multiply the system \eqref{eq:openloop_system} by $P^{-1}$. Defining $\Tilde{B} = P^{-1}B$ and $\Tilde{w}_t = P^{-1}w_t$, the system dynamics become
\begin{equation*}
    \Tilde{x}_{t+1} = \Tilde{A}\Tilde{x}_t + \Tilde{B}u_t +
    \Tilde{w}_t, 
\end{equation*}
which is in the same form as \eqref{eq:openloop_system} but with the system matrix in real Jordan normal form.

For the purposes of controlling this system, it suffices to consider each of the Jordan blocks of $A$ individually. Since the matrix $B$ is invertible, the control of each of these blocks will be identical to the control problem for the full system \eqref{eq:openloop_system}. Therefore, by a slight abuse of notation, we will consider the control of a single mode or Jordan block $A\in\mathbb{R}^{n\times n}$ (which may now be part of a larger system) with the single repeated eigenvalue $\lambda\in\mathbb{C}$. By the real Jordan normal form \cite[Theorem 3.4.1.5]{MatrixAnalysis}, we know then that $A$ takes the form
\begin{equation}\label{eq:realjordan}
\begin{bmatrix}
\lambda & 1 &  & \\
 & \lambda & \ddots & \\
 &  & \ddots & 1 \\
 &  &  & \lambda
\end{bmatrix}
\text{ or }
\begin{bmatrix}
D & I &  & \\
 & D & \ddots & \\
 &  & \ddots & I \\
 &  &  & D
\end{bmatrix}
\end{equation}
where $\lambda\in\mathbb{R}$ and $\lambda\in\cminusr$ respectively. For $\lambda=a+bi\in\cminusr$ above, $I$ is the $2\times 2$ identity matrix and $D$ is the $2\times 2$ matrix
\begin{equation*}
D=
\begin{bmatrix}
a & b\\
-b & a
\end{bmatrix}.
\end{equation*}
We note that for $A$ in either of the forms above, we can describe $\inorm{A}$ quite easily:
\begin{prop}
There are four cases to consider, which are:
\begin{itemize}
    \item If $\lambda\in\mathbb{R}$ and $n=1$ then $\inorm{A} = \abs{\lambda}$.
    \item If $\lambda\in\mathbb{R}$ and $n>1$ then $\inorm{A} = \abs{\lambda} + 1$.
    \item If $\lambda=a+bi\in\cminusr$ and $n=2$ then $\inorm{A} = |a| + |b| \leq \sqrt{2}\abs{\lambda}$.
    \item If $\lambda=a+bi\in\cminusr$ and $n>2$ then $\inorm{A} = |a| + |b| + 1 \leq \sqrt{2}\abs{\lambda} + 1$.
\end{itemize}
\end{prop}
The equalities follow just by observing the rows of $A$ under each assumption. In the complex case, the upper bound follows by Cauchy-Schwarz inequality in $\mathbb{R}^2$. The largest upper bound is $\systemzoomnp$, so in view of \eqref{eq:matrixvector_inequality} we may always write that
\begin{equation}\label{eq:zoomout_rate}
    \inorm{Ax} \leq \systemzoom\inorm{x}.
\end{equation}

Briefly, we remark that the above bound can be improved in the case $n > 1$ by applying an invertible transform $S$. That is, as before we let $\Tilde{A} = S^{-1}AS$ and apply the same ``change-of-view'' transformations $\Tilde{x}_t = Sx_t$, $\Tilde{B} = S^{-1}B$ and $\Tilde{w}_t = S^{-1}w_t$. By left-multiplying the system \eqref{eq:openloop_system} by $S^{-1}$ we arrive at an identical control problem now with the system matrix $\Tilde{A} = S^{-1}AS$.

From here it suffices to determine how small $\inorm{S^{-1}AS}$ can be made over all invertible transformations $S$. At least in the case $\lambda\in\mathbb{R}$ we demonstrate that the infimum is at most $\abslam$ using the following construction.

First, suppose $A$ is a real Jordan block, i.e.\ that of \eqref{eq:realjordan} for $\lambda\in\mathbb{R}$. For any $\eps > 0$ we let $S_\eps \coloneqq \textrm{diag}(1,\eps,\eps^2,...,\eps^{n-1})$ with inverse $S^{-1}_\eps = \textrm{diag}(1,\eps^{-1}, \eps^{-2},...,\eps^{-(n-1)})$. Then,
\begin{equation*}
    S^{-1}_\eps AS_\eps = \begin{bmatrix}
\lambda & \eps &  & \\
 & \lambda & \ddots & \\
 &  & \ddots & \eps \\
 &  &  & \lambda
\end{bmatrix}.
\end{equation*}
Therefore, $\inorm{S^{-1}_\eps AS_\eps} = \abslam + \eps$. By taking $\eps\to 0$ this value is arbitrarily close to $\abslam$. If one allows for complex vector and matrix entries, an identical argument can be made in the case $\lambda\in\mathbb{C}\setminus\mathbb{R}$.

The minimization of $\inorm{A}$ will be relevant to our scheme in terms of the minimum capacity required for stabilization.

\begin{remark}
    It is important to distinguish between the cases $\abslam < 1$ and $\abslam \geq 1$. In the former case, the system is open-loop-stable so stability is easy to achieve. Here, if one employs only the fixed quantization stage of the scheme we present in the next section, then stability of the process is nearly automatic and the asymptotic optimality will follow by simple iterated expectation arguments that are almost identical to what we will present shortly. Therefore, we will ignore the stable case and suppose moving forward that $\abslam \geq 1$.
\end{remark}

\subsection{\SchemeP\ in the Vector Case}\label{sec:joint_schemes}
In this section we describe the joint coding and control scheme of Theorem \ref{thm:ultimate_result}. The scheme presented is very similar to that of the scalar case and is in fact a direct generalization of the joint scheme of Theorem \ref{thm:ultimate_scalar}.

Let $K\geq 2$ be an even integer and suppose $\cbrackets{\Delta_t}_{t=0}^\infty$ is a sequence of positive ``bin sizes'' such that $\Delta_{t+1}$ is a function of only $\Delta_t$ and the indicator random variable $\indicator{\inorm{x_t}\leq \frac{K}{2}\Delta_t}$. Also assume that both the encoder and decoder (controller) know $\Delta_0$. Then so long as the type I quantization $\zoomquantizer(x_t)$ is sent over the channel, it is possible to synchronize knowledge of $\Delta_t$ between the quantizer and the controller since $\inorm{x_t}\leq\frac{K}{2}\Delta_t$ if and only if $\zoomquantizer(x_t) \neq 0$.

The coding scheme is as follows. For $\cbrackets{\Delta_t}_{t=0}^\infty$ as above, we calculate the adaptive quantizer output $\zoomquantizer(x_t)$ and the adaptive system error $e_t \coloneqq x_t - \zoomquantizer(x_t)$. Then for integer $N\geq 2$ we use a fixed type II quantizer $\qfixed$ with bin size $\DeltaN$ as in Section \ref{sec:vecquant} to calculate the fixed quantizer output $\qfixed(e_t)$. We then send $\zoomquantizer(x_t)$ and $\qfixed(e_t)$ across the noiseless channel where the channel capacity is at least,
\begin{equation}\label{eq:rate_vector}
    C = \log_2\parens{K^n+1} + \log_2\parens{\parens{N+1}^n}.
\end{equation}
We estimate the state $x_t$ as $\hat{x}_t \coloneqq \zoomquantizer(x_t) + \qfixed(e_t)$. To mirror the fully observed case, the controller applies the control
\begin{equation*}
    u_t = -B^{-1}A\parens{\zoomquantizer(x_t) + \qfixed(e_t)} = -B^{-1}A\hat{x}_t.
\end{equation*}
The scheme is essentially exactly as illustrated in Figure \ref{fig:two_stage_scheme_scalar} from the scalar case, where $-\frac{a}{b}$ is generalized to $-B^{-1}A$.

The controlled system dynamics resulting from this scheme are
\begin{align}
    x_{t+1} &= A(x_t - \hat{x}_t) + w_t\nonumber\\
    &= A(e_t - \qfixed(e_t)) + w_t.\label{eq:controlled_system}
\end{align}

The update dynamics for $\cbrackets{\Delta_t}_{t=0}^\infty$ are nearly identical to the scalar case. We assume that $K>\inorm{A}$ and choose scalars $\frac{\inorm{A}}{K} < \alpha < 1$, $\rho > \inorm{A}$ and $L>0$. We assume again that $\rho \geq K\alpha$. Choose $\Delta_0 \geq L$ arbitrarily, then for $t\geq1$ the bin update is
\begin{equation}\label{eq:binupdate}
    \Delta_{t+1} = 
    \begin{cases}
        \rho\Delta_t,& \text{if } \inorm{x_t} > \tfrac{K}{2}\Delta_t\\
        \alpha\Delta_t,      & \text{if } \inorm{x_t}\leq \tfrac{K}{2}\Delta_t,\Delta_t\geq L\\
        \Delta_t,           & \text{if } \inorm{x_t}\leq \tfrac{K}{2}\Delta_t,\Delta_t < L.
    \end{cases}
\end{equation}

\begin{prop}
With dynamics \eqref{eq:controlled_system} and \eqref{eq:binupdate}, the process $\ourchain$ is a time-homogeneous Markov chain.
\end{prop}

The motivation for this scheme is as in the scalar case. The adaptive type I quantizer $\zoomquantizer$ will lead to stability in the sense of positive Harris recurrence, and the fixed type II quantizer $\qfixed$ will lead to order-optimal convergence of the system quadratic form $x^{\top}Qx$ under invariant measure as the fixed quantization rate $N$ grows large.

As in the scalar case, we impose Condition \ref{cond:countability} so that the state space for $\Delta_t$ is countable. The state space for $\Delta_t$ is
\begin{equation*}
    \Omega_\Delta \coloneqq \cbrackets{\alpha^j\rho^k\Delta_0 : j,k \in \mathbb{Z}_{\geq 0}}.
\end{equation*}
and the state space for the Markov chain $\ourchain$ is $\mathbb{R}^n\times\Omega_\Delta$. What remains is to specify additional constraints which complete our proposed scheme.

We assume that Condition \ref{cond:beta} holds for some $\beta > 2$. For any $\eps\in(0,\beta-2)$ we finish our construction of \SchemeP\ by requiring that $\rho > (\inorm{A})^{\frac{\beta}{\eps}}$ and specifying the dependence of $\DeltaN$ on $N$ as $\DeltaN = 2N^{-1 + \frac{2}{\beta-\eps}}$ (again, the constant $2$ here is arbitrary).

We impose the following restriction on the minimum bin size which
generalizes Condition~\ref{cond:minbinsize_scalar} to the vector case.

\begin{condition}\label{cond:minbinsize_vector}
The minimum adaptive bin size is at least,
\begin{equation*}
    \alpha L > \frac{\inorm{A}}{K\alpha - \inorm{A}}\DeltaN.
\end{equation*}
\end{condition}

Finally, as in the scalar case, as $C\to\infty$ we keep $K$ fixed and let $N\to\infty$ to take advantage of fixed quantization results at high rates, to be presented shortly.

\subsection{Proof Program for Stability and Convergence}\label{sec:proof_program_vector}
In this section, we outline the proof program for Theorem \ref{thm:ultimate_result}, i.e.\ that \SchemeP\ achieves the exponential rate function $r(C) = 2^{\parens{\frac{4}{\beta-\eps}}\frac{1}{n}C}$. Many of the proofs of intermediate results are tedious and largely mechanical, so we have relegated these to the Appendix.

We first give an intermediate result on high-rate quantizer distortion and then present a high-level proof program of Theorem \ref{thm:ultimate_result} with similar key arguments as those in Section \ref{sec:proof_program_scalar}.

Let $\cbrackets{X_N}_{N=2}^\infty$ be a sequence of random vectors on $\Rn$ and consider the quantizer $\qfixed$ described in Section \ref{sec:vecquant}. We define for $N\geq 2$ the error vectors 
\[Y_N \coloneqq X_N - \qfixed(X_N).\]
\begin{lemma}\label{lem:MSE_1}
Suppose that
\begin{equation*}
    \sup_{N\geq 2}\expec[]{\inorm{X_N}^m} \eqqcolon B_m < \infty
\end{equation*}
for some $m>2$ (not necessarily integer). Set the bin size for the type II quantizer $\qfixed$ as $\DeltaN = 2N^{-1+\frac{2}{m}}$. Then for any positive semidefinite matrix $V\in\mathbb{R}^{n\times n}$ we have
\begin{equation*}
    \trp{V\expec[]{Y_N Y_N^{\top}}} = \landau{N}{N^{-2+\frac{4}{m}}}.
\end{equation*}

\end{lemma}
The proof is mostly mechanical and can be found in the Appendix.

We denote the ``in-view'' set $\Lambda$\footnote{In \cite{YukTAC2010}, this was referred to as the {\it perfectly zoomed} phase set} as,
\begin{equation*}
    \Lambda \coloneqq \cbrackets{(x,\Delta)\in\Rn\times\Omega_\Delta : \inorm{x}\leq\tfrac{K}{2}\Delta},
\end{equation*}
i.e., the set of states $\statep{}$ such that the adaptive quantizer $Q_K^\Delta(x)$ is non-zero.

We will show in Lemma \ref{lem:returntime_tailbound} that states beginning in this set return to it very quickly, in that for constants $h>0$ and $\xi>1$ we have
\begin{equation*}
    \prob[x,\Delta]{\tau_\Lambda \geq k+1} \leq k\Two\parens{ \frac{\Delta}{2}\parens{\frac{h\xi^k-1}{k} -\frac{\DeltaN}{\alpha L}} }
\end{equation*}
which decays very fast in integer $k\geq 1$. Essentially all of the following stability results will follow from repeated application of this inequality and Condition \ref{cond:beta} on the noise process.

\begin{theorem}\label{thm:PHR}
  Under \SchemeP\ with $K > \inorm{A}$, $\ourchain$ is positive Harris
  recurrent for every even $N\geq 2$. Therefore, for every even $N\geq 2$, \SchemeP\ yields a
  unique invariant measure $\piN$ for the process $\ourchain$.
\end{theorem}
\begin{proof}[Sketch of Proof]\let\qed\relax
  We establish $\varphi$-irreducibility and aperiodicity for
  $\ourchain$, where $\varphi$ is the product of the Lebesgue and
  discrete measures on $\Rn\times\Omega_\Delta$. The logarithmic
  function $V(x,\Delta) = c\log_\alpha\Delta$ is shown to satisfy
  Condition \ref{cond:lyapunov_random_base} with $d(x,\Delta)$
  constant and $f\equiv 1$ (i.e., in the form required by Lemma
  \ref{lem:drift_phr}), leading to positive Harris recurrence. A
  complete proof is provided in the Appendix.
\end{proof}

We denote $(\xpi,\Deltapi)\sim\piN$ as the state under invariant measure. This will also induce an invariant measure for the system adaptive error $e_t$, which we denote by $\epi\sim\piN^{\textrm{err}}$.

We have the following ergodicity result, similar to the scalar case.

\begin{prop}\label{prop:technical_ergodicity}
The infinite-horizon second moment and the invariant second moment agree, that is
\begin{equation*}
    \lim_{T\to\infty}\frac{1}{T}\expec[]{\sum_{t=0}^{T-1}x_t^{\top}Qx_t} = \expec[]{\xpiparens^{\top}Q\xpiparens}.
\end{equation*}
\end{prop}
\begin{proof}[Sketch of Proof]
  We are able to show using the drift conditions of Section
  \ref{sec:lyapunov} that functions $g(x,\Delta)$ which are bounded by
  $\inorm{x}^{\beta-\eps}$ satisfy the above ergodicity condition. The
  quadratic form $x^{\top}Qx$ is of the order $\inorm{x}^2$ and hence
  bounded by $\inorm{x}^{\beta-\eps}$, since $\eps < \beta - 2$. This
  will establish the result. A complete proof is provided in the
  Appendix.
\end{proof}

The following proposition is crucial to our proof program.

\begin{prop}\label{prop:totalexpectation}
We have the following characterization of the invariant second moment.
\begin{equation*}
    \expec[]{\xpiparens^{\top}Q\xpiparens} - \trp{Q\Sigma} = \trp{A^{\top}QA\cdot\expec[]{(\epi - \qfixed(\epi))(\epi - \qfixed(\epi))^{\top}}}.
\end{equation*}
\end{prop}
\begin{proof}
The proof is by iterated expectations. Suppose that $\statep{0}\sim\piN$. Let $e_0 = x_0-Q_K^{\Delta_0}(x_0)$ and $x_1 = A(e_0 - \qfixed(e_0)) + Z$, where $Z\sim\eta$ (recall $\eta$ is the distribution of $w_t$). Since we have applied the one-step transition kernel and $\piN$ is the invariant measure, the marginal distributions of $x_0$ and $x_1$ will be identical.

For brevity we denote $s_0 \coloneqq e_0 - \qfixed(e_0)$ so that $x_1 = As_0 + Z$. Supposing that $\expec[]{(\xpi)^{\top}Q(\xpi)} < \infty$ (this follows from system moment results that will be stated shortly), we then have by invariance and iterated expectations that
\begin{align*}
    \expec[]{x_0^{\top}Qx_0} &= \expec[]{x_1^{\top}Qx_1} = \expec[]{ \expec[]{x_1^{\top}Qx_1\ \big|\ s_0} }\\
    &= \expec[]{ \expec[]{(As_0 + Z)^{\top}Q(As_0 + Z)\ \big|\ s_0} }\\
    &= \expec[]{ s_0^{\top}A^{\top}QAs_0 + 2\expec[]{Z}^{\top}QAs_0 + Z^{\top}QZ }\\
    &= \expec[]{s_0^{\top}A^{\top}QAs_0} + \expec[]{Z^{\top}QZ}\\
    &= \trp{A^{\top}QA\cdot \expec[]{s_0s_0^{\top}}} + \trp{Q\expec[]{ZZ^{\top}}}\\
    &= \trp{A^{\top}QA\cdot \expec[]{(e_0 - \qfixed(e_0))(e_0 - \qfixed(e_0))^{\top}}} + \trp{Q\Sigma}.
\end{align*}
Rearranging the above equality completes the proof. Note that above we used the property that $Z\sim\eta$ is zero-mean.
\end{proof}

\begin{remark}
For the type of scheme we present here (that is, using two-stage adaptive and fixed uniform quantization), if one is able to show that
\begin{equation*}
    \trp{A^{\top}QA\cdot\expec[]{(\epi - \qfixed(\epi))(\epi - \qfixed(\epi))^{\top}}} = \landau{N}{\frac{R(N)}{N^2}}
\end{equation*}
for some function $R(N)$ then it follows from the above two propositions (and the fact that $N$ is a linear function of $2^{\frac{1}{n}C}$, recalling \eqref{eq:rate_vector} and that $K$ is fixed) that the scheme achieves the rate function
\begin{equation*}
    r(C) = R\parens{ (K^n+1)^{-\frac{1}{n}}2^{\frac{1}{n}C} - 1 }.
\end{equation*}
That is, we have
\begin{equation}\label{eq:NtoC}
    \lim_{T\to\infty}\frac{1}{T}\expec[]{\sum_{t=0}^{T-1}x_t^{\top}Qx_t} - \trp{Q\Sigma} = \landau{C}{\frac{R\parens{ (K^n+1)^{-\frac{1}{n}}2^{\frac{1}{n}C} - 1 }}{2^{\frac{2}{n}C}}}.
\end{equation}
Depending on the specific function $R(N)$ in question, this expression can be simplified (this will be the case in our analysis).
\end{remark}
We have the following uniform stability result under invariant
measure.

\begin{lemma}\label{lem:schemeAmoment}
Under \SchemeP, the invariant system error $\epi$ has finite $(\beta-\eps)$-th moment uniformly in $N\geq 2$. That is,
\begin{equation}\label{eq:superrormoments_vector}
    \sup_{N\geq 2}\expec[]{\inorm{\epi}^{\beta-\eps}} < \infty.
\end{equation}
\end{lemma}
\begin{proof}[Sketch of Proof]
With Lyapunov functions $V(x,\Delta)$ and $d(x,\Delta)$ proportional to $\Delta^{\beta-\eps}$, an appropriate ``in-view'' set $C$ and constant $b$ we show that these functions satisfy the drift condition \eqref{eq:driftcriteria_foraset}. In particular, we do this for $f$ proportional to $\inorm{x}^{\beta-\eps}$ and $f$ proportional to $\Delta^{\beta-\eps}$ which, with the Lyapunov parameters independent of $N$ leads to results of the form \eqref{eq:superrormoments_vector} by Lemma \ref{lem:drift_moment} for the invariant state and adaptive bin size. A simple invariance argument finishes the result. A detailed proof is provided in the Appendix.
\end{proof}

We note that $A^{\top}QA$ is positive semidefinite (by positive definiteness of $Q$). This allows us to use Lemma \ref{lem:MSE_1}. Finally, we prove our ultimate result.

\begin{proof}[Proof of Theorem \ref{thm:ultimate_result}]
Lemma \ref{lem:schemeAmoment} allows us to use Lemma \ref{lem:MSE_1} with the sequence of random vectors $\cbrackets{\epi}_{N=2}^\infty$ from which we obtain,
\begin{equation*}
    \trp{A^{\top}QA\cdot\expec[]{(\epi - \qfixed(\epi))(\epi - \qfixed(\epi))^{\top}}} = \landau{N}{\frac{ N^{\frac{4}{\beta-\eps}} }{ N^2 }}
\end{equation*}
and so by the earlier remark with $R(N) = N^{\frac{4}{\beta-\eps}}$, \SchemeP\ achieves the rate function
\begin{equation*}
    r(C) = \parens{ (K^n+1)^{-\frac{1}{n}} 2^{\frac{1}{n}C} - 1 }^{\frac{4}{\beta - \eps}} = \landau{C}{2^{\parens{\frac{4}{\beta-\eps}}\frac{1}{n}C}}.
\end{equation*}
This proves Theorem \ref{thm:ultimate_result} and completes our proof program.
\end{proof}

%\iffalse

\section{Simulation Results}\label{sec:simulation}
In this section, we provide an example simulation to illustrate our results. Consider the linear scalar system,
\begin{equation*}
    x_{t+1} = 1.2x_t + u_t + w_t.
\end{equation*}
Here the noise process $\cbrackets{w_t}_{t=0}^\infty$ is i.i.d.\
with the marginal distribution of $w_0 = 4Z$,  where $Z$ is
``Bucklew-Gallagher'' distributed with pdf \eqref{eq:bg_pdf} with 
$\delta = 2$. This marginal distribution is related to the typical
Pareto distribution in that
$\frac{1}{4}|w_t| + 1 \sim \textrm{Pareto}(1, 4)$.

We remark that $w_t$ admits finite moments $\beta$ only of order $\beta < 4$, whereas all moments $\beta \geq 4$ are infinite. These random variables are thus badly behaved in the sense of having heavy tails. Despite this, it is clear that $\cbrackets{w_t}_{t=0}^\infty$ satisfies Condition \ref{cond:beta_scalar} for all $\beta\in(2,4)$. We are interested in the minimization of
\begin{equation*}
    \limsup_{T\to\infty}\frac{1}{T}\expec[]{\sum_{t=0}^{T-1}x_t^2}
\end{equation*}
across a noiseless discrete channel of capacity $C$ . This asymptotic
performance has a lower bound of $\expec[]{w_0^2} = \frac{16}{3}$.

Let $K=2$ (the least even integer such that $K > |a| = 1.2$). Then for
even $N\geq 2$, the capacity of our channel is
\begin{equation*}
    C = \log_2(K+1) + \log_2(N+1) = \log_2 3 + \log_2(N+1).
\end{equation*}
Let $\alpha = \frac{3}{4}$, $\rho = \parens{\frac{4}{3}}^3$ and $L=9$ be the adaptive quantization parameters. Note that $\alpha^3\rho = 1$ so that Condition \ref{cond:countability} is satisfied. Let $\beta = 3.95$ and $\eps = 0.95$. Finally, we let $\DeltaN = 2N^{-1 + \frac{2}{\beta - \eps}} = 2N^{-\frac{1}{3}}$ be the bin size for the fixed quantizer $\qfixed$ for all even $N\geq 2$.

Note that $\rho > (1.2)^{79/19} = |a|^{\frac{\beta}{\eps}}$ and that $\rho \geq K\alpha = \frac{3}{2}$. Also note that Condition \ref{cond:minbinsize_scalar} is satisfied for all $N\geq 2$ since $\alpha L = \frac{27}{4} > \frac{8}{2^{1/3}} = \frac{|a|}{K\alpha - |a|}\Delta_{(2)}$.

Therefore, with $\beta - \eps = 3$, by employing \SchemeP\ we are guaranteed by Theorem \ref{thm:ultimate_scalar} that,
\begin{equation*}
    \limsup_{T\to\infty}\frac{1}{T}\expec[]{\sum_{t=0}^{T-1}x_t^2} - \frac{16}{3} = \landau{C}{2^{\parens{-2 + \frac{4}{\beta-\eps}}C}} = \landau{C}{2^{-\frac{2}{3}C}}.
\end{equation*}
Equivalently, in the sense of Definition \ref{def:convergence_scalar}
the scheme presented here achieves the exponential rate function $r(C)
= 2^{\frac{4}{3}C}$.

\begin{figure}[h]
    \centering
    \includegraphics[width=0.8\linewidth]{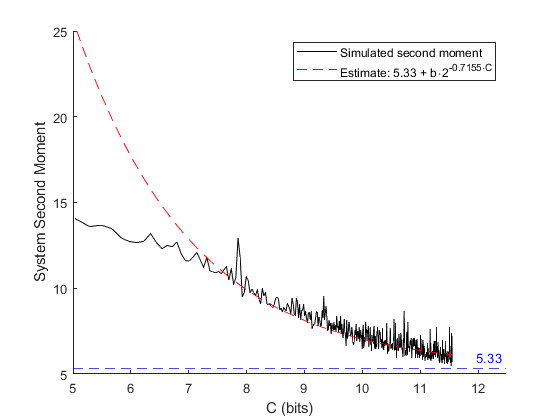}
    \caption{Convergence to the optimum $\expec[]{w_0^2} = \frac{16}{3}$. The order of convergence is approximately $\landau{C}{2^{-0.7155C}}$, which is just slightly better than the expected convergence $\landau{C}{2^{-\frac{2}{3}C}}$. The constant $b$ is on the order of $2^8$.}
    \label{fig:simulation}
\end{figure}

\begin{figure}[h]
    \centering
    \includegraphics[width=0.8\linewidth]{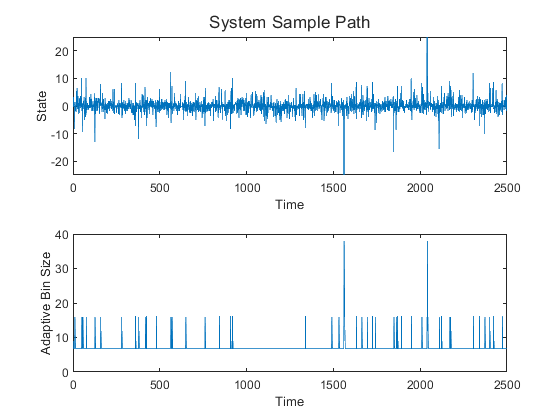}
    \caption{A sample path for the system with $N=100$ fixed quantizer bins. The system state and adaptive bin size are displayed over time.}
    \label{fig:sample_path}
\end{figure}

Using this scheme with the parameters above, the system was run for
all $N\in\cbrackets{10, 12, ...\ 1000}$ and the average second moment
was recorded. We remark that the settling time for the average second
moment was observed to vary dramatically between small and large
values of $N$. Thus, instead of a fixed simulation length $T$ for each trial of $N$, a variable stopping time was employed. Let
\begin{equation*}
	S_T \coloneqq \frac{1}{T} \sum_{t=0}^{T-1}x_t^2
\end{equation*}
be the average second moment up to time $T>0$. A small convergence threshold $\epsilon>0$ and an integer settling time $T^*>0$ are defined, and the simulation stops when for $T^*$ consecutive time stages we have $|S_{T+1} - S_T| < \epsilon$. That is, the average second moment varies no more than $\epsilon$ across a single time stage for $T^*$ consecutive time stages. For the simulations presented here, the convergence threshold used is $\epsilon = 10^{-4}$ and the settling time is $T^* = 10^4$.

The second moment achieved in each trial (i.e., for each $N$) is plotted against the capacity $C = \log_2 3 + \log_2(N+1)$ below in Figure \ref{fig:simulation}.

% \begin{figure}[h]
%     \centering
%     \includegraphics[width=0.8\linewidth]{moments.png}
%     \caption{Convergence to the optimum $\expec[]{w_0^2} = \frac{16}{3}$. The order of convergence is approximately $\landau{C}{2^{-0.7155C}}$, which is just slightly better than the expected convergence $\landau{C}{2^{-\frac{2}{3}C}}$. The constant $b$ is on the order of $2^8$.}
%     \label{fig:simulation}
% \end{figure}

% \begin{figure}[h]
%     \centering
%     \includegraphics[width=0.8\linewidth]{sample_path.png}
%     \caption{A sample path for the system with $N=100$ fixed quantizer bins. The system state and adaptive bin size are displayed over time.}
%     \label{fig:sample_path}
% \end{figure}

Since we expect that the optimality gap will converge to zero at a rate like $2^{-\delta C}$ for some $\delta > 0$, it is reasonable to expect that the logarithm of the optimality gap is approximately linear in $C$ with slope $-\delta$. Therefore, an estimate of the order of convergence was obtained by performing a linear regression between $C$ and the logarithm of the optimality gap. This estimate on the order of convergence can be seen in Figure \ref{fig:simulation}. Finally, a finite segment of a sample path for the Markov chain $\ourchain$ is displayed in Figure \ref{fig:sample_path} for the case $N=100$.

Finally, we remark on the convergence rate to the invariant measure $\piN$. In \cite[Theorem IV.2]{RYLTAC2015} it is shown that in the scalar case with Gaussian noise, using only the adaptive part of the coding scheme, the Markov chain $\ourchain$ is ``geometrically ergodic'' in that for every $(x,\Delta)\in\mathbb{R}\times\Omega_\Delta$, the $n$-stage transition kernels $P^n((x,\Delta),\cdot)$ converge to $\piN(\cdot)$ under the total variation metric at a rate $r^n$ for some $r<1$.

Much of the analysis needs to be re-derived in the more general case
we consider here, but one can still follow the primary line of
argument in the proof of \cite[Theorem IV.2]{RYLTAC2015}. In
particular, invoking \cite[Theorem III.7]{RYLTAC2015} with the
Lyapunov function $V(x,\Delta) = \Delta^2$ and using our
super-geometric bound on the tail probabilities
$\prob[x,\Delta]{\tau_\Lambda \geq k+1}$, Lemma
\ref{lem:returntime_tailbound}, one can show  that the scheme presented
here also ensures that the Markov chain $\ourchain$ is geometrically
ergodic for each $N\geq 2$ (in this proof, which we have not worked
out in detail,  a  careful argument is needed to ensure that the
geometric rate is uniformly bounded over all $N\geq 2$). This
difficulty is of similar flavor that we  resolved by a careful choice of Lyapunov
parameters in Lemma \ref{lem:schemeAmoment}).

%\fi

\section{Conclusion and Future Work}
In conclusion, we have constructed joint coding and
control schemes for networked control systems of the form
\eqref{eq:openloop_system} which are asymptotically optimal in the
sense that, as
the data rate grows without bound, the system second moment converges
to the classical optimum with an explicit rate of convergence. The
techniques in this paper build on prior work in this context, in
particular by the use of random-time Lyapunov drift conditions to
establish key stability results.

There are several potential directions for future work. It would be
interesting to explore optimality of schemes for \emph{non-linear}
systems. The two-stage scheme approach seems applicable here, supposing that one can establish stability and
ergodicity results by the adaptive part of the code. For instance, in
\cite[Theorem 5.1]{yuksel2015stability}, a nonlinear system with
``sub-linear'' dynamics and additive Gaussian noise is shown to admit
stability when subjected to an adaptive zooming quantization scheme
similar to what is presented in \cite{YukTAC2010,YukMeynTAC2010} and
the adaptive stage of the scheme presented here. It seems feasible
that employing a two-stage (adaptive and fixed) uniform quantization
scheme for systems of this type would, by some extra analysis, lead to
convergence results similar to what we have presented here.

One direction is the relaxation of the invertibility assumption on the
control matrix $B$. In general, it is possible to achieve stability in
the sense of positive Harris recurrence and \emph{finite} system
moments with just the assumption that the pair $(A,B)$ is controllable
\cite[Theorem 2.2]{AndrewJohnstonReport}. Controllability is a natural
relaxation of the invertibility assumption to pursue for the kinds of
linear systems considered here. Nonetheless, for optimality arguments
on rates of decay in the distortion as the rate increases, further
analysis is needed.

A further possibility is to consider stricter conditions on the noise
process than Condition \ref{cond:beta}. For instance, one might assume
that the noise tails are dominated by exponential decay (i.e., the
noise has sub-exponential distribution), and seek to construct schemes
which, in the sense of Definition \ref{def:convergence_vector},
achieve rate functions that are much better than exponential (e.g.,
polynomial in $C$). The two-part coding scheme presented here is well
suited to generalizations of this type, because the bound of Lemma
\ref{lem:returntime_tailbound} holds quite generally (e.g., with few
assumptions on the noise process) and repeated use of this bound and
Condition \ref{cond:beta} leads to most of our key results. Thus,
using this bound with a stricter condition on the noise process may
lead to stronger stability and optimality results.

\section*{Appendix: Proofs}

\begin{proof}[Proof of Proposition \ref{prop:fully_observed_optimal}]
We show that in the fully observed setup, the optimal control policy which minimizes \eqref{eq:ergodic_cost} is $u_t = -B^{-1}Ax_t$, achieving an optimal cost of $\trp{Q\Sigma}$, where $\trp{\cdot}$ is the trace operator. Let $v_t \coloneqq Ax_t + Bu_t$ so that $x_{t+1} = v_t + w_t$, then under any policy $\gamma$ we have
\begin{align}
    \limsup_{T\to\infty}\frac{1}{T}\expec[]{\sum_{t=0}^{T-1}x_t^{\top}Qx_t} &= \limsup_{T\to\infty}\frac{1}{T}\expec[]{\sum_{t=0}^{T-1} x_{t+1}^{\top}Qx_{t+1}}\nonumber\\
    &= \limsup_{T\to\infty}\frac{1}{T}\expec[]{\sum_{t=0}^{T-1}(v_t+w_t)^{\top}Q(v_t+w_t)}\nonumber\\
    &= \limsup_{T\to\infty}\frac{1}{T}\expec[]{\sum_{t=0}^{T-1}v_t^{\top}Qv_t + w_t^{\top}Qw_t}\nonumber\\
    &\geq \limsup_{T\to\infty}\frac{1}{T}\expec[]{\sum_{t=0}^{T-1}w_t^{\top}Qw_t}\label{eqstep:fo_optimal_vt}\\
    &= \expec[]{w_0^{\top}Qw_0} = \trp{Q\expec[]{w_0w_0^{\top}}} = \trp{Q\Sigma}\nonumber
\end{align}
by positive definiteness of $Q$ and that $w_t$ is i.i.d.\ zero-mean. In the case $u_t = -B^{-1}Ax_t$ we have $v_t = 0$ and so \eqref{eqstep:fo_optimal_vt} is an equality, establishing optimality.
\end{proof}

\begin{proof}[Proof of Lemma \ref{lem:ultimaterate_scalar}]
Suppose that $b=1$ without loss of generality so that we are considering control of the following system,
\begin{equation*}
    x_{t+1} = ax_t + u_t + w_t
\end{equation*}
across a discrete noiseless channel of capacity $C$ bits using an arbitrary joint coding and control scheme $\phi_C$. First, note we may assume that for all $C$ sufficiently large,
\begin{equation}\label{eqstep:state_limsup_finite}
    \limsup_{T\to\infty}\frac{1}{T}\expec[]{\sum_{t=0}^{T-1}x_t^2} < \infty.
\end{equation}
If this fails to be true for the joint scheme $\phi_C$ then \eqref{eq:scalar_capacity_bound} will trivially hold. Note that \eqref{eqstep:state_limsup_finite} implies that
\begin{equation}\label{eqstep:state_liminf_finite}
    \liminf_{T\to\infty}\expec[]{x_T^2} < \infty.
\end{equation}

We will follow the core arguments in the proof of \cite[Theorem 11.3.2]{YukselBasarBook} with the key difference of not assuming that $\lim_{t\to\infty} \expec[]{x_t^2}$ exists.

Denote the received channel output at time $t\geq 0$ as $q'_t$. Let $D_t \coloneqq \expec[]{(x_t+\frac{1}{a}u_t)^2}$, $d_t \coloneqq \expec[]{x_t^2}$ and note that the system update gives us $D_t = \frac{1}{a^2}(d_{t+1} - \sigma^2)$.

We note that the discrete noiseless channel we consider is a special (noise-free discrete memoryless) case of ``Class A'' channels \cite[Definition 8.5.1]{YukselBasarBook}. Such channels have their capacity characterized by limits of the directed mutual information between the channel input and output. In particular, following the beginning of the proof of \cite[Theorem 8.5.2]{YukselBasarBook}, we find that our channel capacity $C$ satisfies
\begin{equation}\label{eq:capacity_characterization}
    C \geq \limsup_{T\to\infty} \frac{1}{T}\sum_{t=0}^{T-1} I(x_t;q'_t|q'_{[0,t-1]}),
\end{equation}
where $I(x;y|z)$ is the conditional mutual information. Following the analysis in the proof of \cite[Theorem 11.3.2]{YukselBasarBook} (essentially exactly that on page 392) and invoking a conditional version of the entropy-power inequality \cite[Lemma 5.3.2]{YukselBasarBook} we are able to arrive at the following bound,
\begin{align*}
    \frac{1}{T}\sum_{t=0}^{T-1} I(x_t;q'_t|q'_{[0,t-1]}) &\geq \frac{1}{T}\sum_{t=0}^{T-1}\frac{1}{2}\log_2\parens{a^2 + \frac{\sigma^2}{D_t}} + \frac{1}{T}\parens{h(x_0) - h(ax_{T-1}+w_{T-1}\ |\ q'_{[0,T-1]})}\\
    &\geq \frac{1}{2}\log_2\parens{a^2 + \frac{\sigma^2}{\frac{1}{T}\sum_{t=0}^{T-1}D_t}} + \frac{1}{T}\parens{h(x_0) - h(ax_{T-1}+w_{T-1}\ |\ q'_{[0,T-1]})}
\end{align*}
by the convexity of $x\mapsto\log(1+\frac{1}{x})$. Above, $h(\cdot)$ and $h(\cdot\ |\ \cdot)$ are the regular and conditional differential entropies, respectively. Therefore, we arrive at the following bound on the channel capacity in view of \eqref{eq:capacity_characterization},
\begin{align*}
    C &\geq \limsup_{T\to\infty} \parens{ \frac{1}{2}\log_2\parens{ a^2 + \frac{\sigma^2}{\frac{1}{T}\sum_{t=0}^{T-1}D_t} } + \frac{1}{T}\parens{ h(x_0) - h(ax_{T-1} + w_{T-1}\ |\ q'_{[0,T-1]}) } }\\
    &\geq \liminf_{T\to\infty} \frac{1}{2}\log_2\parens{ a^2 + \frac{\sigma^2}{\frac{1}{T}\sum_{t=0}^{T-1}D_t} } + \limsup_{T\to\infty}\frac{1}{T}\parens{ h(x_0) - h(ax_{T-1} + w_{T-1}\ |\ q'_{[0,T-1]}) }.
\end{align*}
We will consider the two limits above separately. Let
\begin{equation*}
    d \coloneqq \limsup_{T\to\infty} \frac{1}{T}\sum_{t=0}^{T-1}d_t = \limsup_{T\to\infty}\frac{1}{T}\expec[]{\sum_{t=0}^{T-1}x_t^2} < \infty
\end{equation*}
and since $D_t = \frac{1}{a^2}(d_{t+1} - \sigma^2)$ we have that $\limsup_{T\to\infty}\frac{1}{T}\sum_{t=0}^{T-1}D_t = \frac{1}{a^2}(d-\sigma^2)$.

Then since $\log\parens{1+\frac{1}{x}}$ is continuous and monotone decreasing in $x>0$ we have
\begin{align*}
    \liminf_{T\to\infty} \frac{1}{2}\log_2\parens{ a^2 + \frac{\sigma^2}{\frac{1}{T}\sum_{t=0}^{T-1}D_t} } &= \frac{1}{2}\log_2\parens{ a^2 + \frac{\sigma^2}{\limsup_{T\to\infty}\frac{1}{T}\sum_{t=0}^{T-1}D_t} }\\
    &= \frac{1}{2}\log_2\parens{a^2 + \frac{a^2\sigma^2}{d-\sigma^2}}\\
    &= \frac{1}{2}\log_2\parens{\frac{a^2 d}{d-\sigma^2}}.
\end{align*}
Therefore, we have that
\begin{equation}
    C \geq \frac{1}{2}\log_2\parens{\frac{a^2 d}{d-\sigma^2}} + \limsup_{T\to\infty}\frac{1}{T}\parens{ h(x_0) - h(ax_{T-1} + w_{T-1}\ |\ q'_{[0,T-1]}) }
\end{equation}
and what remains is to bound this limit supremum. We then have,
\begin{align}
    \limsup_{T\to\infty}\frac{1}{T}&\parens{ h(x_0) - h(ax_{T-1} + w_{T-1}\ |\ q'_{[0,T-1]}) }\nonumber\\
    &= -\liminf_{T\to\infty} \frac{1}{T}h(ax_{T-1} + w_{T-1}\ |\ q'_{[0,T-1]})\nonumber\\
    &= -\liminf_{T\to\infty} \frac{1}{T}h(x_T - bu_{T-1}\ |\ q'_{[0,T-1]})\nonumber\\
    &= -\liminf_{T\to\infty}\frac{1}{T}h(x_T\ |\ q'_{[0,T-1]})\label{eqstep:cond_diff_invariant}\\
    &\geq -\liminf_{T\to\infty}\frac{1}{T}h(x_T).\label{eqstep:cond_diff_bound}
\end{align}
Above, \eqref{eqstep:cond_diff_invariant} follows since $u_t$ is constant given $q'_{[0,t]}$ and differential entropy is translation invariant. \eqref{eqstep:cond_diff_bound} follows since conditioning reduces differential entropy.

We now show that $\liminf_{T\to\infty}\frac{1}{T}h(x_T) \leq 0$.

Note that since the noise process is added to the state at every time stage and has a pdf $\eta$ which is positive everywhere on $\mathbb{R}$, the state $x_t$ will also have a positive-everywhere density. Furthermore, the state $x_t$ will also have finite variance $\expec[]{x_t^2} - \expec[]{x_t}^2$ at each time stage.

It is known that for a distribution over $\mathbb{R}$ with specified variance $S$, the Gaussian distribution maximizes differential entropy at $\frac{1}{2}\log_2\parens{2\pi eS}$. Therefore we have,
\begin{align*}
    \liminf_{T\to\infty}\frac{1}{T}h(x_T) &\leq \liminf_{T\to\infty}\frac{1}{2T}\log_2\parens{2\pi e\parens{ \expec[]{x_T^2} - \expec[]{x_T}^2 }}\\
    &\leq \liminf_{T\to\infty}\frac{1}{2T}\log_2\parens{2\pi e \expec[]{x_T^2}}\\
    &\leq \liminf_{T\to\infty}\frac{1}{T}\log_2\parens{\expec[]{x_T^2}} = 0,
\end{align*}
where the final limit is zero by \eqref{eqstep:state_liminf_finite}.

Therefore, we have the following ultimate bound on the channel capacity,
\begin{equation*}
    C\geq\frac{1}{2}\log_2\parens{\frac{a^2d}{d-\sigma^2}}.
\end{equation*}
Recall that $d = \limsup_{T\to\infty}\frac{1}{T}\expec[]{\sum_{t=0}^{T-1}x_t^2}$. Rearranging the above inequality for the optimality gap $d-\sigma^2$ we find that
\begin{equation*}
    \limsup_{T\to\infty}\frac{1}{T}\expec[]{\sum_{t=0}^{T-1}x_t^2} -
    \sigma^2 \geq \frac{a^2\sigma^2}{2^{2C} - a^2}, 
\end{equation*}
which completes the proof.
\end{proof}

\begin{proof}[Proof of Lemma \ref{lem:MSE_1}]
Define the event $F$ by
\begin{equation*}
    F = \cbrackets{\pnorm{X_N}{\infty} \leq \tfrac{1}{2}N\DeltaN}.
\end{equation*}
To start, the $(i,j)$-th component of $Y_NY_N^{\top}$ is $Y_N^iY_N^j$. Therefore,
\begin{align}
    \abs{\sbrackets{\expec[]{Y_NY_N^{\top}}}_{ij}} &= \abs{\expec[]{Y_N^iY_N^j}} = \abs{ \expec[]{Y_N^iY_N^j\one_F} + \expec[]{Y_N^iY_N^j\one_{F^C}} }\nonumber\\
    &\leq \abs{\expec[]{Y_N^iY_N^j\one_F}} + \abs{\expec[]{Y_N^iY_N^j\one_{F^C}}}\nonumber\\
    &\leq \expec[]{\abs{Y_N^iY_N^j}\one_F} + \expec[]{\abs{Y_N^iY_N^j}\one_{F^C}}\label{eqstep:mse_twocases}
\end{align}
by triangle inequality and Jensen's inequality. Since $F$ is such that $\abs{X_N^k} \leq \frac{1}{2}N\DeltaN$ for all $1\leq k \leq n$, it follows by construction of $\qfixed$ that $\abs{Y_N^k} \leq \frac{1}{2}\DeltaN$. Therefore,
\begin{equation}\label{eqstep:mse_case1}
    \expec[]{\abs{Y_N^iY_N^j}\one_F} \leq \expec[]{\parens{\tfrac{1}{2}\DeltaN}^2\one_F} \leq \parens{\tfrac{1}{2}\DeltaN}^2 = \tfrac{1}{4}\DeltaN^2.
\end{equation}
Note that $\abs{Y_N^k} \leq \abs{X_N^k}$, then for the second expectation of \eqref{eqstep:mse_twocases} we have by H\"older's inequality that
\begin{align}
    \expec[]{\abs{Y_N^iY_N^j}\one_{F^C}} &\leq \expec[]{\abs{X_N^iX_N^j}\one_{F^C}} \leq \expec[]{\abs{X_N^iX_N^j}^{\frac{m}{2}}}^{\frac{2}{m}}\expec[]{\one_{F^C}}^{1-\frac{2}{m}}\nonumber\\
    &= \expec[]{\abs{X_N^iX_N^j}^{\frac{m}{2}}}^{\frac{2}{m}}\prob[]{\inorm{X_N} > \tfrac{1}{2}N\DeltaN}^{\frac{m-2}{m}}.\label{eqstep:mse_holder1}
\end{align}
Recall that we suppose $\sup_{N\geq 2} \expec[]{\inorm{X_N}^m} \eqqcolon B_m < \infty$. It follows by Cauchy-Schwarz inequality that the expectation in \eqref{eqstep:mse_holder1} is bounded as
\begin{align}
    \expec[]{\abs{X_N^iX_N^j}^{\frac{m}{2}}}^{\frac{2}{m}} &\leq \expec[]{\abs{X_N^i}^m}^{\frac{1}{m}}\expec[]{\abs{X_N^j}^m}^{\frac{1}{m}}\nonumber\\
    &\leq \expec[]{\inorm{X_N}^m}^{\frac{2}{m}} \leq (B_m)^{\frac{2}{m}}.\label{eqstep:mse_cs1}
\end{align}
Note that by Markov's inequality, we have for $u>0$ that
\begin{align*}
    \prob[]{\inorm{X_N} > u} &= \prob[]{\inorm{X_N}^m > u^m} \leq \expec[]{\inorm{X_N}^m}u^{-m} \leq B_mu^{-m}
\end{align*}
and so since $\frac{1}{2}N\DeltaN = N^{\frac{2}{m}}$, the tail probability in \eqref{eqstep:mse_holder1} can be bounded as,
\begin{align}
    \prob[]{\inorm{X_N} > \tfrac{1}{2}N\DeltaN}^{\frac{m-2}{m}} &\leq \parens{B_m\parens{N^{\frac{2}{m}}}^{-m}}^{\frac{m-2}{m}}\nonumber\\
    &= (B_m)^{\frac{m-2}{m}}N^{-2 + \frac{4}{m}}.\nonumber
\end{align}
Combining this with \eqref{eqstep:mse_holder1} and \eqref{eqstep:mse_cs1} we find that
\begin{align*}
    \expec[]{\abs{Y_N^iY_N^j}\indicator{F^C}} \leq (B_m)^{\frac{2}{m} + \frac{m-2}{m}}N^{-2+\frac{4}{m}} = B_m N^{-2+\frac{4}{m}} = \landau{N}{\DeltaN^2},
\end{align*}
which, combined with \eqref{eqstep:mse_case1} yields that
\begin{equation*}
    \abs{\sbrackets{\expec[]{Y_NY_N^{\top}}}_{ij}} = \landau{N}{\DeltaN^2}.
\end{equation*}
Now that we have demonstrated that each component of $\expec[]{Y_NY_N^{\top}}$ satisfies the desired bound, the proof completes as follows. For any matrix $V\in\mathbb{R}^{n\times n}$ we have by triangle inequality and Jensen's inequality that
\begin{align}
    \abs{\trp{V\expec[]{Y_NY_N^{\top}}}} &= \abs{\sum_{i=1}^n\sum_{j=1}^nV_{ij}\sbrackets{\expec[]{Y_NY_N^{\top}}}_{ji}} = \abs{\sum_{i=1}^n\sum_{j=1}^nV_{ij}\expec[]{Y_N^iY_N^j}}\nonumber\\
    &\leq \sum_{i=1}^n\sum_{j=1}^n \big|V_{ij}\big| \abs{\expec[]{Y_N^iY_N^j}}\nonumber\\
    &\leq \parens{\max_{k,l}\abs{V_{kl}}}\sum_{i=1}^n\sum_{j=1}^n\abs{\expec[]{Y_N^iY_N^j}} = \landau{N}{\DeltaN^2}.\label{eqstep:mse_lemma_final}
\end{align}
The proof can be completed by noting that for any two positive semidefinite matrices $P,Q$ we have $\trp{PQ} \geq 0$. Note that for any random vector $X$, $\expec[]{XX^{\top}}$ is positive semidefinite. Therefore, for any positive semidefinite matrix $V$ it follows that
\begin{equation*}
    \trp{V\expec[]{Y_NY_N^{\top}}} = \abs{ \trp{V\expec[]{Y_NY_N^{\top}}} }
\end{equation*}
and so the proof concludes in view of \eqref{eqstep:mse_lemma_final}.
\end{proof}

We now proceed with proving stability of $\ourchain$ in the sense of both positive Harris recurrence as well as moment stability. To show positive Harris recurrence, we will need to demonstrate that our chain is irreducible and that an appropriate class of sets are small. The proofs of these two technical results (Propositions \ref{prop:aperiodicity} and \ref{prop:smallsets}) are tedious and so we only provide proof sketches for brevity. For details,  see \cite{keelerThesis}.

\begin{prop}\label{prop:aperiodicity}
The process $\ourchain$ is $\varphi$-irreducible and aperiodic, where $\varphi$ is the product of the Lebesgue and discrete measures on $\Rn\times\Omega_\Delta$.
\end{prop}
\begin{proof}[Sketch of Proof]
The condition of aperiodicity is strictly stronger than that of irreducibility, so it suffices to show only aperiodicity. Loosely, from any initial state $\statep{0}$ and a target set $B$ with $\varphi(B) > 0$ there exists some $n_0 > 0$ such that for every $n\geq n_0$, we can drive the system from $\statep{0}$ to $B$ in $n$ time stages with positive probability.

The reason this is possible is that we may ``hold'' the bin size $\Delta_t$ constant when it is less than $L$, and so by careful management we may direct the bin size $\Delta_n$ to anywhere in the state space $\Omega_\Delta$ we like (with a big enough $n_0$). In particular, we can direct the system to $B$ so long as the state $x_n$ falls into a subset of positive Lebesgue measure at time $n$. That this happens with positive probability is because the state is convolved with noise at every time stage that has a positive-everywhere pdf $\eta$.
\end{proof}

We denote the in-view set as
\begin{equation}
    \Lambda \coloneqq \cbrackets{(x,\Delta)\in\Rn\times\Omega_\Delta : \inorm{x}\leq\tfrac{K}{2}\Delta}
\end{equation}
and in light of Section \ref{sec:lyapunov} we will develop results concerning the return time $\tau_\Lambda$.

\begin{prop}\label{prop:smallsets}
For the Markov chain $\ourchain$, bounded subsets of $\Lambda$ are small.
\end{prop}
\begin{proof}[Sketch of Proof]
    We can show that subsets of $\Lambda$ containing only one bin size are $1$-small. This follows mainly since such sets are bounded and $\Delta_1$ is known deterministically, so the next-stage pdf is bounded from below by a sub-probability measure in terms of the initial set and the noise pdf $\eta$.

    It follows by aperiodicity that any finite union of small sets is small, and since bounded subsets of $\Lambda$ are the finite union of subsets with one bin size (which are $1$-small), these sets are small (though generally not $1$-small).
\end{proof}

The following proposition will prove to be remarkably useful for the remainder of our proof program.

\begin{prop}\label{prop:iid_infinity_bound}
Let $\cbrackets{z_t}_{t=0}^\infty$ be an i.i.d.\ sequence of nonnegative random variables. For any $b>0$ and integer $k\geq 1$ we have
\begin{equation*}
    \prob[]{\sum_{t=0}^{k-1}z_t > b} \leq k\prob[]{z_0 > \tfrac{b}{k}}.
\end{equation*}
\end{prop}
\begin{proof}
Since $\cbrackets{z_t}_{t=0}^\infty$ is identically distributed, we have for $k\geq 1$ that
\begin{align*}
    \prob[]{\sum_{t=0}^{k-1}z_t > b} &\leq \prob[]{\bigcup_{t=0}^{k-1}\cbrackets{z_t > \tfrac{b}{k}}} \leq \sum_{t=0}^{k-1}\prob[]{z_t > \tfrac{b}{k}} = k\prob[]{z_0 > \tfrac{b}{k}}.
\end{align*}
\end{proof}

\begin{corollary}\label{cor:iid_moment_bound}
Let $\cbrackets{z_t}_{t=0}^\infty$ be an i.i.d.\ sequence of nonnegative random variables. Then for any (real) $m>0$ and integer $k\geq 1$ we have
\begin{equation*}
    \expec[]{\parens{ \sum_{t=0}^{k-1}z_t}^m} \leq k^{m+1}\expec[]{z_0^m}.
\end{equation*}
\end{corollary}
\begin{proof}
Using the tail formula for expectation of a nonnegative random variable and the previous proposition, we have
\begin{align*}
    \expec[]{\parens{ \sum_{t=0}^{k-1}z_t}^m} &= \int_0^\infty \prob[]{\parens{ \sum_{t=0}^{k-1}z_t}^m > u}du = \int_0^\infty \prob[]{\sum_{t=0}^{k-1}z_t > u^{\frac{1}{m}}}du\\
    &\leq k\int_0^\infty \prob[]{z_0 > \frac{u^{\frac{1}{m}}}{k}} = k\int_0^\infty \prob[]{(kz_0)^m > u}du\\
    &= k\expec[]{(kz_0)^m} = k^{m+1}\expec[]{z_0^m}.
\end{align*}
\end{proof}

\begin{lemma}\label{lem:returntime_tailbound}
Define the constants
\begin{equation*}
    \xi \coloneqq \frac{\rho}{\inorm{A}},\quad h \coloneqq \frac{K\alpha}{\rho}.
\end{equation*}
For $\statep{}\in\Lambda$ we have for any $k\geq 1$ that
\begin{align}
    \prob[\state{}]{\tau_\Lambda \geq k+1} &\leq k\Two\parens{\frac{\Delta}{2}\parens{\frac{h\xi^k-1}{k} - \frac{\DeltaN}{\alpha L}}}.\label{eq:pzt_tail_bound2}
\end{align}
\end{lemma}
\begin{proof}
Starting from $\statep{}\in\Lambda$, let $e_0 \coloneqq x - Q_K^{\Delta}(x)$. We define the map $S:\Rn\to\Rn$ by
\begin{equation*}
    S(u) \coloneqq u - \qfixed(u).
\end{equation*}
Note that $\inorm{S(u)} \leq \inorm{u} + \frac{\DeltaN}{2}$ for all $u\in\Rn$. The correction term $\frac{\DeltaN}{2}$ accounts for the possibility that $\inorm{u} < \frac{\DeltaN}{2}$.

We construct the following "zoom-out" process $y_t$. Let $y_0=e_0$, and for $t\geq 0$ let
\begin{equation}\label{eq:zoomout_process}
    y_{t+1} = AS(y_t) + w_t.
\end{equation}
Then the following holds for all $1\leq t \leq \tau_\Lambda$:
\begin{equation*}
    y_t = x_t,\quad\Delta_t = \rho^{t-1}\alpha\Delta.
\end{equation*}
Therefore, it follows that for $k\geq 1$,
\begin{align*}
    \cbrackets{\tau_\Lambda \geq k+1} &\ =\ \bigcap_{t=1}^k\cbrackets{\statep{t}\not\in\Lambda} \ =\ \bigcap_{t=1}^k\cbrackets{\inorm{x_t} > \tfrac{K}{2}\Delta_t}\\
    &\ =\ \bigcap_{t=1}^k\cbrackets{\inorm{y_t} > \tfrac{\Delta}{2}K\alpha\rho^{t-1}}\\
    &\ \subseteq\  \cbrackets{\inorm{y_k} > \tfrac{\Delta}{2}K\alpha\rho^{k-1}} \ =\ \cbrackets{\inorm{y_k} > \tfrac{\Delta}{2}h\rho^k}.
\end{align*}
We may further relax this event as
\begin{align}
    \cbrackets{\inorm{y_k} > \tfrac{\Delta}{2}h\rho^k} &\ =\ \cbrackets{\inorm{AS(y_{k-1}) + w_{k-1}} > \tfrac{\Delta}{2}h\rho^k}\nonumber\\
    &\ \subseteq\ \cbrackets{\inorm{AS(y_{k-1})} + \inorm{w_{k-1}} > \tfrac{\Delta}{2}h\rho^k}\label{eqstep:inorm_triangle}\\
    &\ \subseteq\ \cbrackets{\inorm{A}\inorm{S(y_{k-1})} + \inorm{w_{k-1}} > \tfrac{\Delta}{2}h\rho^k}\label{eqstep:zoombound_jordan}\\
    &\ =\ \cbrackets{\inorm{S(y_{k-1})} + \inorm{A}^{-1}\inorm{w_{k-1}} > \tfrac{\Delta}{2}h\xi\rho^{k-1}}\nonumber\\
    &\ \subseteq\ \cbrackets{\inorm{S(y_{k-1})} + \inorm{w_{k-1}} > \tfrac{\Delta}{2}h\xi\rho^{k-1}}\label{eqstep:systemzoon_morethanone}\\
    &\ \subseteq\ \cbrackets{\inorm{y_{k-1}} + \inorm{w_{k-1}} + \tfrac{\DeltaN}{2} > \tfrac{\Delta}{2}h\xi\rho^{k-1}}.\label{eqstep:S_contractive}
\end{align}
Above, \eqref{eqstep:inorm_triangle} follows by triangle inequality, \eqref{eqstep:zoombound_jordan} holds in light of \eqref{eq:zoomout_rate}, \eqref{eqstep:systemzoon_morethanone} holds since $\inorm{A} \geq \abslam \geq 1$ and \eqref{eqstep:S_contractive} holds since $\inorm{S(u)}\leq\inorm{u} + \frac{\DeltaN}{2}$.

The steps \eqref{eqstep:inorm_triangle} through \eqref{eqstep:systemzoon_morethanone} can be repeated $k-1$ more times to ultimately obtain that the event \eqref{eqstep:systemzoon_morethanone} implies
\begin{equation}\label{eqstep:inviewtime_bound_almostdone}
    \cbrackets{\inorm{y_0} + \sum_{t=0}^{k-1}\parens{\inorm{w_t}+\tfrac{\DeltaN}{2}} > \tfrac{\Delta}{2}h\xi^k}.
\end{equation}
Using the fact that $y_0 = e_0$ and that starting in-view, $\inorm{e_0} \leq \frac{\Delta}{2}$, we find that
\begin{align*}
    \cbrackets{\tau_\Lambda \geq k+1} \ \subseteq\ \cbrackets{\sum_{t=0}^{k-1}\parens{\inorm{w_t}+\tfrac{\DeltaN}{2}} > \tfrac{\Delta}{2}\parens{h\xi^k-1}}.
\end{align*}
Therefore,
\begin{equation*}
    \prob[\state{}]{\tau_\Lambda\geq k+1} \leq \prob[\state{}]{\sum_{t=0}^{k-1}\parens{\inorm{w_t}+\tfrac{\DeltaN}{2}} > \tfrac{\Delta}{2}\parens{h\xi^k-1}}.
\end{equation*}
Note that since $\alpha > \frac{\inorm{A}}{K}$ we have $h\xi^k > 1$ for $k\geq 1$. Therefore, $\tfrac{\Delta}{2}\parens{h\xi^k-1} > 0$ so we may apply Proposition \ref{prop:iid_infinity_bound} to the sequence of i.i.d. random variables $\cbrackets{\inorm{w_t}+\tfrac{\DeltaN}{2}}_{t=0}^\infty$ to find that
\begin{align*}
    \prob[\state{}]{\tau_\Lambda\geq k+1} &\leq k\prob[x,\Delta]{\inorm{w_0} + \tfrac{\DeltaN}{2} > \frac{\tfrac{\Delta}{2}\cdot(h\xi^k-1)}{k}}\\
    &= k\prob[x,\Delta]{\inorm{w_0} > \frac{\Delta}{2}\cdot\frac{h\xi^k-1}{k} - \frac{\DeltaN}{2}}\\
    &= k\Two\parens{\frac{\Delta}{2}\cdot\frac{h\xi^k-1}{k} - \frac{\DeltaN}{2}}.
\end{align*}
Finally, note that with $\Delta \geq \alpha L$ we have,
\begin{align*}
    \frac{\Delta}{2}\cdot\frac{h\xi^k-1}{k} - \frac{\DeltaN}{2} &= \frac{\Delta}{2}\parens{\frac{h\xi^k-1}{k} - \frac{\DeltaN}{\Delta}} \geq \frac{\Delta}{2}\parens{\frac{h\xi^k-1}{k} - \frac{\DeltaN}{\alpha L}},
\end{align*}
which, since $\Two\parens{\cdot}$ is nonincreasing, proves the claimed bound \eqref{eq:pzt_tail_bound2}.

Briefly, we justify that the argument to $\Two\parens{\cdot}$ in \eqref{eq:pzt_tail_bound2} is \emph{strictly positive} for all $k\geq 1$. First, note that Condition \ref{cond:minbinsize_vector} is equivalent to,
\begin{equation}\label{eq:pzt_bin_bound}
    \frac{\DeltaN}{\alpha L} < h\xi-1.
\end{equation}

Since $\alpha > \frac{\inorm{A}}{K}$ and $\rho \geq K\alpha$ we have $h\in(\frac{1}{\xi}, 1]$. It is relatively straightforward (if tedious) to verify that for $\xi > 1$ and $h\in(\frac{1}{\xi}, 1]$ that the function
\begin{equation*}
    s\mapsto \frac{h\xi^s-1}{s}
\end{equation*}
is monotone increasing for real $s>0$ (e.g. by careful analysis of its derivative). In particular, it follows that $\frac{h\xi^k-1}{k} \geq h\xi-1$ for integer $k\geq 1$. Therefore, it follows from \eqref{eq:pzt_bin_bound} (i.e., from Condition \ref{cond:minbinsize_vector}) that
\begin{align*}
    \frac{\Delta}{2}\parens{\frac{h\xi^k-1}{k} - \frac{\DeltaN}{\alpha L}} > \frac{\Delta}{2}\parens{\frac{h\xi^k-1}{k} - (h\xi-1)} \geq 0.
\end{align*}
\end{proof}

\begin{corollary}\label{cor:return_inview_bothschemes}
For $\statep{}\in\Lambda$, the in-view return time satisfies
\begin{align*}
    &\sup_{(x,\Delta)\in\Lambda}\expec[x,\Delta]{\tau_\Lambda} < \infty\\
    &\lim_{\Delta\to\infty} \expec[x,\Delta]{\tau_\Lambda - 1} = 0.
\end{align*}
\end{corollary}
\begin{proof}
Note that by Condition \ref{cond:beta} , we have $\expec[]{\inorm{w_0}^2} < \infty$ (this follows since $\beta > 2$). Therefore by Markov's inequality we have for $u>0$ that
\begin{equation*}
    \Two(u) \leq \expec[]{\inorm{w_0}^2}u^{-2}.
\end{equation*}
We use the tail formula for the expectation of a nonnegative discrete random variable and Lemma \ref{lem:returntime_tailbound} to then find that for $\statep{}\in\Lambda$,
\begin{align*}
    \expec[\state{}]{\tau_\Lambda} &= \sum_{k=1}^\infty \prob[\state{}]{\tau_\Lambda \geq k} = 1 + \sum_{k=1}^\infty \prob[\state{}]{\tau_\Lambda \geq k+1}\\
    &\leq 1 + \sum_{k=1}^\infty k\Two\parens{\frac{\Delta}{2}\parens{\frac{h\xi^k-1}{k} - \frac{\DeltaN}{\alpha L}}}\\
    &\leq 1 + \sum_{k=1}^\infty k\expec[]{\inorm{w_0}^2}\parens{\frac{\Delta}{2}\parens{\frac{h\xi^k-1}{k} - \frac{\DeltaN}{\alpha L}}}^{-2}\\
    &= 1 + \Delta^{-2}\cdot\sum_{k=1}^\infty 4\expec[]{\inorm{w_0}^2}k^3\parens{h\xi^k-1-\frac{\DeltaN}{\alpha L}k}^{-2},
\end{align*}
where the series converges because each term is $\landau{k}{k^3\xi^{-2k}}$ for $\xi>1$. We remark that $\Delta \geq \alpha L$ implies that the above is uniformly bounded (when one replaces $\Delta$ by $\alpha L$) and that as $\Delta\to\infty$ the above is $1 + \landau{\Delta}{\Delta^{-2}}$, which proves both claims.
\end{proof}

Finally, we prove positive Harris recurrence.

\begin{proof}[Proof of Theorem \ref{thm:PHR}]
Following the remark after Condition \ref{cond:lyapunov_random_base}, we will satisfy the drift condition \eqref{eq:driftcriteria_foraset} for Lemma \ref{lem:drift_phr}. In this case, we set the following.
\begin{align*}
    &V(\state{}) = c\log_\rho\Delta,\\
    &d(\state{}) \equiv d = \sup_{\statep{}\in\Lambda}\expec[\state{}]{\tau_\Lambda},\\
    &f(\state{}) \equiv 1,\\
    &C = \Lambda\cap\cbrackets{\statep{}\in\Rn\times\Omega_\Delta : \Delta\leq D},\\
    &b = c(d-1) + c\log_\rho\alpha + d,
\end{align*}
where $d$ is finite by Corollary \ref{cor:return_inview_bothschemes}, $c\geq \frac{2d}{\log_\rho(\frac{1}{\alpha})}$ and $D$ is such that $\Delta > D$ implies that $\expec[\state{}]{\tau_\Lambda - 1} \leq \tfrac{1}{2}\log_\rho(\frac{1}{\alpha})$ (such a $D$ exists by Corollary \ref{cor:return_inview_bothschemes}). We note that $C$ is a small set by Proposition \ref{prop:smallsets}.

First, notice that by construction we have for any $\statep{}\in\Lambda$ that $\expec[\state{}]{\tau_\Lambda} \leq d$, so the second inequality of \eqref{eq:driftcriteria_foraset} is satisfied. We now show that the first inequality also holds, which in this case is that for any $\statep{}\in\Lambda$,
\begin{equation}\label{eq:phr_tosatisfy}
    \expec[\state{}]{V(\state{\tau_\Lambda})} - V(\state{}) \leq - d + b\indicator{\statep{}\in C}.
\end{equation}
First, note that for $1\leq t\leq \tau_\Lambda$ we have that $\Delta_t = \rho^{t-1}\alpha\Delta$, so,
\begin{align*}
    V(\state{\tau_\Lambda}) &= c\log_\rho\Delta_{\tau_\Lambda} = c\log_\rho\parens{\rho^{\tau_\Lambda-1}\alpha\Delta}\\
    &= c(\tau_\Lambda-1) + c\log_\rho\alpha + c\log_\rho\Delta\\
    &= c(\tau_\Lambda-1) + c\log_\rho\alpha + V(\state{}).
\end{align*}
Thus we have that
\begin{equation}\label{eqstep:phr_ts1}
    \expec[\state{}]{V(\state{\tau_\Lambda})} - V(\state{}) = c\expec[\state{}]{\tau_\Lambda - 1} + c\log_\rho\alpha.
\end{equation}
First, suppose $\Delta>D$ (so that $\statep{}\not\in C$), then by construction of $D$ and $c$, \eqref{eqstep:phr_ts1} becomes
\begin{align*}
    c\expec[\state{}]{\tau_\Lambda-1} + c\log_\rho\alpha &\leq \tfrac{1}{2}c\log_\rho(\tfrac{1}{\alpha}) + c\log_\rho\alpha = -\tfrac{1}{2}c\log_\rho(\tfrac{1}{\alpha})\\
    &\leq -\frac{1}{2}\parens{\frac{2d}{\log_\rho(\frac{1}{\alpha})}}\log_\rho(\tfrac{1}{\alpha})\\
    &= -d = -d + b\indicator{\statep{}\in C},
\end{align*}
which satisfies \eqref{eq:phr_tosatisfy}. Next, suppose that $\Delta \leq D$, then by construction of $b$ and $d$, \eqref{eqstep:phr_ts1} becomes
\begin{align*}
    c\expec[\state{}]{\tau_\Lambda - 1} + c\log_\rho\alpha &\leq c(d-1) + c\log_\rho\alpha\\
    &= -d + b = -d + b\indicator{\statep{}\in C}
\end{align*}
and so \eqref{eq:phr_tosatisfy} is satisfied over $\Lambda$. It follows then by Lemma \ref{lem:drift_phr} that $\ourchain$ is positive Harris recurrent, provided that we can show that $\prob[\state{}]{\tau_C < \infty} = 1$ for all $\statep{}\in \Rn\times\Omega_\Delta$. We do this to complete the proof.

Since we have proven the drift condition \eqref{eq:driftcriteria_foraset} holds at all $\statep{}\in\Lambda$, it follows from the proof of Lemma \ref{lem:drift_phr} that for all $\statep{}\in\Lambda$ we have
\begin{equation*}
    \expec[\state{}]{\tau_C} \leq V(\state{}) + b < \infty
\end{equation*}
so we must have that $\prob[\state{}]{\tau_C < \infty} = 1$ for all $\statep{}\in\Lambda$. From here, it suffices to verify that all states return to $\Lambda$ in finite time. This is automatic for initial $\statep{}\in\Lambda$ by Corollary \ref{cor:return_inview_bothschemes}. Now suppose that $\statep{}\in\Lambda^C$. By nearly identical arguments to those of Lemma \ref{lem:returntime_tailbound}, we can show that
\begin{align*}
    \prob[\state{}]{\tau_\Lambda \geq k+1} &\leq \prob[\state{}]{\sum_{t=0}^{k-1}\parens{\inorm{w_t} + \tfrac{\DeltaN}{2}} > \tfrac{\Delta}{2}K\xi^k - \inorm{x}}.
\end{align*}
For $k$ sufficiently large we have $\tfrac{\Delta}{2}K\xi^k > \inorm{x}$, so we may apply Proposition \ref{prop:iid_infinity_bound} which yields that for $k$ sufficiently large,
\begin{align*}
    \prob[\state{}]{\tau_\Lambda \geq k+1} &\leq k\Two\parens{\frac{\tfrac{\Delta}{2}K\xi^k - \inorm{x}}{k} - \frac{\DeltaN}{2}},
\end{align*}
which converges to zero as $k$ grows large ($\Two(u) = \landau{u}{u^{-2}}$, as in the proof of Corollary \ref{cor:return_inview_bothschemes}). Therefore, we must have $\prob[\state{}]{\tau_\Lambda = \infty} = 0$, which completes the proof.
\end{proof}

\begin{proof}[Proof of Proposition \ref{prop:technical_ergodicity}]
First, we note that if $\lambda^*$ is the maximum eigenvalue of $Q$ then the following holds for arbitrary $x\in\Rn$,
\begin{equation}\label{eq:quadraticformbound}
    x^{\top}Qx \leq \lambda^* \pnorm{x}{2}^2 \leq
    \lambda^*\parens{\sqrt{n}\inorm{x}}^2 = n\lambda^*\inorm{x}^2, 
\end{equation}
where the first bound follows by properties of positive semidefinite matrices and the second follows by \eqref{eq:pnorm_inequality}. Since $\eps < \beta - 2$ it follows that $x^{\top}Qx$ is bounded by $\inorm{x}^{\beta-\eps}$ in the sense of Lemma \ref{lem:drift_moment}.

We will show in the proof of Proposition \ref{prop:state_moment_sup} that Condition \ref{cond:lyapunov_random_base} holds in the form required by Lemma \ref{lem:drift_moment} with the function $f(x,\Delta) = C\inorm{x}^{\beta-\eps}$, for a constant $C>0$. Since $x^{\top}Qx$ is bounded by $\inorm{x}^{\beta-\eps}$, the proof completes in view of the ergodicity result \eqref{eq:driftmomentergodic}.
\end{proof}

To complete the proof program, we must show the moment condition of Lemma \ref{lem:schemeAmoment} holds. The rest of the Appendix is dedicated to this proof. We first have some intermediate results.

\begin{prop}\label{prop:returntime_schemeA}
Under \SchemeP, for $\statep{}\in\Lambda$ the return time $\tau_\Lambda$ satisfies the following independently of $N$.
\begin{align*}
    &\sup_{\statep{}\in\Lambda} \expec[]{\parens{\rho^{\beta-\eps}}^{\tau_\Lambda-1}} < \infty\\
    &\lim_{\Delta\to\infty} \expec[]{\parens{\rho^{\beta-\eps}}^{\tau_\Lambda-1}} = 1.
\end{align*}
\end{prop}
\begin{proof}
Note that by Condition \ref{cond:beta} with $\beta > 2$ and Markov's inequality we have for $u>0$ that
\begin{equation*}
    \Two(u) \leq \expec[]{\inorm{w_0}^\beta}u^{-\beta}.
\end{equation*}
Let $r \coloneqq \rho^{\beta - \eps}$ for brevity. Then we have by Lemma \ref{lem:returntime_tailbound} that
\begin{align*}
    \expec[\state{}]{r^{\tau_\Lambda - 1}} &= \sum_{k=0}^\infty \prob[\state{}]{\tau_\Lambda = k+1}r^k \leq 1 + \sum_{k=1}^\infty \prob[\state{}]{\tau_\Lambda \geq k+1}r^k\\
    &\leq 1 + \sum_{k=1}^\infty k\Two\parens{\frac{\Delta}{2}\parens{\frac{h\xi^k-1}{k} - \frac{\DeltaN}{\alpha L}}}r^k\\
    &\leq 1 + \sum_{k=1}^\infty kr^k\expec[]{\inorm{w_0}^\beta}\parens{\frac{\Delta}{2}\parens{\frac{h\xi^k-1}{k} - \frac{\DeltaN}{\alpha L}}}^{-\beta}\\
    &= 1 + \Delta^{-\beta}\cdot\sum_{k=1}^\infty 2^\beta\expec[]{\inorm{w_0}^\beta}k^{\beta+1}r^k\parens{h\xi^k-1-\frac{\DeltaN}{\alpha L}k}^{-\beta}.
\end{align*}
Above, the series converges since the summand is $\landau{k}{k^{\beta+1}\parens{r\xi^{-\beta}}^k}$, where $r\xi^{-\beta} < 1$ (this is ensured by the assumption that $\rho > (\inorm{A})^{\frac{\beta}{\eps}})$. Then, we have shown that
\begin{equation*}
    \expec[\state{}]{r^{\tau_\Lambda - 1}} = 1 + \landau{\Delta}{\Delta^{-\beta}},
\end{equation*}
which (in light of the fact that $\Delta \geq \alpha L$) completes the proof.
\end{proof}

\begin{prop}\label{prop:state_moment_sup}
Under \SchemeP, the invariant state $\xpi$ has finite $(\beta-\eps)$-th moment uniformly in $N\geq 2$. That is,
\begin{equation*}
    \sup_{N\geq 2} \expec[]{\inorm{\xpi}^{\beta-\eps}} < \infty.
\end{equation*}
\end{prop}
\begin{prop}\label{prop:bin_moment_sup}
Under \SchemeP, the invariant adaptive bin size $\Deltapi$ has finite $(\beta-\eps)$-th moment uniformly in $N\geq 2$. That is,
\begin{equation*}
    \sup_{N\geq 2} \expec[]{(\Deltapi)^{\beta-\eps}} < \infty.
\end{equation*}
\end{prop}
We will prove these two results shortly via Lyapunov drift arguments, but using them we will first provide a short proof of Lemma \ref{lem:schemeAmoment}.

\begin{proof}[Proof of Lemma \ref{lem:schemeAmoment}]
Let $(\xpi,\Deltapi)\sim\piN$ and let $\epi = \xpi - Q_K^{\Deltapi}(\xpi)$. It follows that $\epi$ is distributed as the invariant system adaptive error. Note that for any $(x,\Delta)$ we have the inequality $\inorm{x - Q_K^\Delta(x)} \leq \inorm{x} + \frac{\Delta}{2}$, where the extra $\frac{\Delta}{2}$ term accounts for the possibility that $\inorm{x} < \frac{\Delta}{2}$. Therefore, we have that
\begin{equation*}
    \inorm{\epi} = \inorm{\xpi - Q_K^{\Deltapi}(\xpi)} \leq \inorm{\xpi} + \frac{\Deltapi}{2}
\end{equation*}
and so,
\begin{align}
    \expec[]{\inorm{\epi}^{\beta-\eps}} &\leq \expec[]{\parens{\inorm{\xpi} + \frac{\Deltapi}{2}}^{\beta-\eps}}\nonumber\\
    &\leq 2^{\beta-\eps}\expec[]{\inorm{\xpi}^{\beta-\eps}} + \expec[]{(\Deltapi)^{\beta-\eps}}.\label{eqstep:error_moments_interms}
\end{align}
Above, the second bound follows from the general fact that for any $\theta > 0$ and nonnegative random variables $X,Y$ we have
\begin{equation}\label{eq:sum_power_bound}
    \expec[]{(X+Y)^\theta} \leq 2^\theta\parens{\expec[]{X^\theta} + \expec[]{Y^\theta}},
\end{equation}
which follows by convexity of $u\mapsto u^\theta$ when $\theta \geq 1$ (via Jensen's inequality) and by sub-additivity of $u\mapsto u^\theta$ when $\theta < 1$.

Therefore we have by \eqref{eqstep:error_moments_interms} that
\begin{align*}
    \sup_{N\geq 2}\expec[]{\inorm{\epi}^{\beta-\eps}} \leq 2^{\beta-\eps}\sup_{N\geq 2}\expec[]{\inorm{\xpi}^{\beta-\eps}} + \sup_{N\geq 2}\expec[]{(\Deltapi)^{\beta-\eps}} < \infty,
\end{align*}
which is finite by Propositions \ref{prop:state_moment_sup} and \ref{prop:bin_moment_sup}.
\end{proof}
All that remains is to prove Propositions \ref{prop:state_moment_sup} and \ref{prop:bin_moment_sup}, which we do here via Lyapunov drift arguments.

\begin{proof}[Proof of Proposition \ref{prop:state_moment_sup}]
Recall that we must show that
\begin{equation*}
    \sup_{N\geq 2} \expec[]{\inorm{\xpi}^{\beta-\eps}} < \infty.
\end{equation*}
We will do this via random-time Lyapunov drift arguments, in particular Lemma \ref{lem:drift_moment} using the drift condition \eqref{eq:driftcriteria_foraset}. Let $r \coloneqq \rho^{\beta-\eps}$ as in Proposition \ref{prop:returntime_schemeA}. We set the following:
\begin{align*}
    &V(\state{}) = \Delta^{\beta-\eps},\\
    &d(\state{}) = s\Delta^{\beta-\eps},\\
    &f(\state{}) = cs\parens{\tfrac{2}{K}\inorm{x}}^{\beta-\eps},\\
    &C = \Lambda\cap\cbrackets{\statep{}\in\Rn\times\Omega_\Delta : \Delta\leq D},\\
    &b = D^{\beta-\eps}\parens{\alpha^{\beta-\eps}\sup_{\state{}\in\Lambda}\expec[\state{}]{r^{\tau_\Lambda-1}}-1+s},
\end{align*}
where $s\in (0, 1-\alpha^{\beta-\eps})$ is arbitrary, $D>0$ is such that
\begin{equation*}
    \statep{}\in\Lambda\textrm{ and }\Delta > D \implies \alpha^{\beta-\eps}\expec[\state{}]{r^{\tau_\Lambda-1}}-1+s \leq 0
\end{equation*}
(such a $D$ exists by Proposition \ref{prop:returntime_schemeA}), $b$ is finite by Proposition \ref{prop:returntime_schemeA}, and $c>0$ is a sufficiently small constant (which will be specified shortly). We will show that the drift condition \eqref{eq:driftcriteria_foraset} holds for the choices above.

First, note that for $\statep{}\in\Lambda$ we have that $\Delta_{\tau_\Lambda} = \rho^{\tau_\Lambda-1}\alpha\Delta$, so
\begin{align}
    \expec[\state{}]{V(\state{\tau_\Lambda})}& - V(\state{}) + d(\state{})\nonumber\\
    &= \expec[\state{}]{\parens{\rho^{\tau_\Lambda-1}\alpha\Delta}^{\beta-\eps}} - \Delta^{\beta-\eps} + s\Delta^{\beta-\eps}\nonumber\\
    &= \Delta^{\beta-\eps}\parens{\alpha^{\beta-\eps}\expec[\state{}]{r^{\tau_\Lambda-1}}-1+s}.\label{eqstep:schemeAmoment1}
\end{align}
For $\Delta\leq D$, by construction of $b$ we have that \eqref{eqstep:schemeAmoment1} is bounded by $b$. For $\Delta > D$, by construction of $D$ we have that \eqref{eqstep:schemeAmoment1} is nonpositive. In either case, the first drift inequality is satisfied.

Next we show that the second inequality of \eqref{eq:driftcriteria_foraset} holds, when $c$ is sufficiently small. To start, we have that for $\statep{}\in\Lambda$ that
\begin{align}
    \expec[\state{}]{\sum_{t=0}^{\tau_\Lambda-1}f(\state{t})} &= f(\state{}) + \expec[\state{}]{\sum_{t=1}^{\tau_\Lambda-1}f(\state{t})}\nonumber\\
    &= cs\parens{\tfrac{2}{K}\inorm{x}}^{\beta-\eps} + cs\parens{\tfrac{2}{K}}^{\beta - \eps}\expec[\state{}]{\sum_{t=1}^{\tau_\Lambda-1}\inorm{x_t}^{\beta-\eps}}\nonumber\\
    &\leq cs\Delta^{\beta-\eps} + cs\parens{\tfrac{2}{K}}^{\beta - \eps}\expec[\state{}]{\sum_{t=1}^{\tau_\Lambda-1}\inorm{x_t}^{\beta-\eps}}\nonumber\\
    &= cd(\state{}) + cs\parens{\tfrac{2}{K}}^{\beta - \eps}\expec[\state{}]{\sum_{t=1}^{\tau_\Lambda-1}\inorm{x_t}^{\beta-\eps}}.\label{eqstep:schemeAmoment2}
\end{align}
The remainder of the proof will be dedicated to showing that there exists a constant $M>0$ so that uniformly over $\statep{}\in\Lambda$ we have,
\begin{align}\label{eqstep:schemeAmomentMain}
    \expec[\state{}]{\sum_{t=1}^{\tau_\Lambda-1}\inorm{x_t}^{\beta-\eps}} \leq M\Delta^{\beta-\eps}
\end{align}
so that we may bound \eqref{eqstep:schemeAmoment2} as
\begin{align*}
    \eqref{eqstep:schemeAmoment2} &\leq cd(\state{}) + cs\parens{\tfrac{2}{K}}^{\beta-\eps}M\Delta^{\beta-\eps} = c\parens{1 + \parens{\tfrac{2}{K}}^{\beta-\eps}M}d(\state{}).
\end{align*}
It then suffices to take $c=(1 + (\tfrac{2}{K})^{\beta-\eps}M)^{-1}$, and the second drift inequality holds. Therefore, we dedicate the remainder of this proof to establishing \eqref{eqstep:schemeAmomentMain}.

As in the proof of Lemma \ref{lem:returntime_tailbound}, the ``zoom-out'' process $y_t$ defined by $y_0 = e_0$, $y_{t+1} = AS(y_t) + w_t$ agrees with $x_t$ for $1\leq t\leq\tau_\Lambda$. Therefore, we have that
\begin{align}
    \expec[\state{}]{\sum_{t=1}^{\tau_\Lambda-1}\inorm{x_t}^{\beta-\eps}} &= \expec[\state{}]{\sum_{t=1}^{\tau_\Lambda-1}\inorm{y_t}^{\beta-\eps}}\nonumber\\
    &= \expec[\state{}]{\sum_{t=1}^\infty\indicator{\tau_\Lambda \geq t+1}\inorm{y_t}^{\beta-\eps}}\nonumber\\
    &= \sum_{t=1}^\infty \expec[\state{}]{\indicator{\tau_\Lambda\geq
      t+1}\inorm{y_t}^{\beta-\eps}} , \label{eqstep:schemeAmoment3}
\end{align}
where the exchange of summation and expectation is justified by the monotone convergence theorem \cite[Theorem 1.26]{BigRudin}.

Let $a \coloneqq \inorm{A}$ for brevity. Now, choose $q\in\parens{\frac{\beta}{\eps},\ \frac{\beta}{\beta-\eps}\log_a\xi}$ arbitrarily (this is a valid interval by the condition that $\rho > a^{\frac{\beta}{\eps}}$). We let $p$ be the H\"older conjugate of $q$ (i.e., $\frac{1}{p} + \frac{1}{q} = 1$) and apply H\"older's inequality to \eqref{eqstep:schemeAmoment3}.
\begin{equation}\label{eqstep:schemeAmoment4}
    \eqref{eqstep:schemeAmoment3} \leq \sum_{t=1}^\infty \prob[\state{}]{\tau_\Lambda\geq t+1}^{\frac{1}{q}}\expec[\state{}]{\inorm{y_t}^{p(\beta-\eps)}}^{\frac{1}{p}}.
\end{equation}
We will consider each factor of the summand separately. First, by Lemma \ref{lem:returntime_tailbound} and familiar arguments we have
\begin{align*}
    \prob[\state{}]{\tau_\Lambda\geq t+1} &\leq t\Two\parens{\frac{\Delta}{2}\parens{\frac{h\xi^t-1}{t} - \frac{\DeltaN}{\alpha L}}}\\
    &\leq t\expec[]{\inorm{w_0}^\beta}\parens{\frac{\Delta}{2}\parens{\frac{h\xi^t-1}{t} - \frac{\DeltaN}{\alpha L}}}^{-\beta}\\
    &\leq \expec[]{\inorm{w_0}^\beta}2^\beta(\alpha L)^{-\beta}\cdot t^{\beta+1}\parens{h\xi^t-1 - \frac{\DeltaN}{\alpha L}t}^{-\beta}
\end{align*}
so that for series convergence,
\begin{equation}\label{eqstep:schemeAseries1}
    \prob[\state{}]{\tau_\Lambda\geq t+1}^{\frac{1}{q}} = \landau{t}{t^{\frac{\beta+1}{q}}\parens{\xi^{-\frac{\beta}{q}}}^t}
\end{equation}
and this term is $\landau{\Delta}{1}$. Next, we consider the second summand factor of \eqref{eqstep:schemeAmoment4}. For brevity, let $m \coloneqq p(\beta-\eps)$. Since $\inorm{Av} \leq a\inorm{v}$ and $\inorm{S(v)} \leq \inorm{v} + \frac{\DeltaN}{2}$, we have for $t\geq 1$ that
\begin{equation*}
    \inorm{y_t} = \inorm{AS(y_{t-1}) + w_{t-1}} \leq a\inorm{y_{t-1}} + \inorm{w_{t-1}} + a\cdot\tfrac{\DeltaN}{2}.
\end{equation*}
Repeating this argument $t-1$ times yields 
\begin{align*}
    \inorm{y_t} &\leq a^t\parens{\inorm{y_0} + \sum_{i=0}^{t-1}a^{-i}\parens{\frac{\inorm{w_i}}{a} + \tfrac{\DeltaN}{2}}}
\end{align*}
and since $a>1$ and $\inorm{y_0} = \inorm{e_0} \leq \frac{\Delta}{2}$ we find that
\begin{equation*}
    \inorm{y_t} \leq a^t\parens{\frac{\Delta}{2} + \sum_{i=0}^{t-1}\parens{\inorm{w_i} + \tfrac{\DeltaN}{2}}}.
\end{equation*}
Therefore, we have by Corollary \ref{cor:iid_moment_bound} and repeated application of \eqref{eq:sum_power_bound} that
\begin{align*}
    \expec[\state{}]{\inorm{y_t}^m} &\leq a^{mt}\expec[\state{}]{\parens{\frac{\Delta}{2} + \sum_{i=0}^{t-1}\parens{\inorm{w_i} + \tfrac{\DeltaN}{2}}}^m}\\
    &\leq a^{mt}\parens{ \Delta^m + 2^m\expec[]{\parens{\sum_{i=0}^{t-1}\parens{\inorm{w_i} + \tfrac{\DeltaN}{2}}}^m} }\\
    &\leq a^{mt}\parens{ \Delta^m + 2^m t^{m+1}\expec[]{\parens{\inorm{w_0} + \tfrac{\DeltaN}{2}}^m} }\\
    &\leq a^{mt}\parens{ \Delta^m + t^{m+1}4^m\parens{\expec[]{\inorm{w_0}^m} + \parens{\tfrac{\DeltaN}{2}}^m} }\\
    &\leq \Delta^m a^{mt}\parens{ 1 + t^{m+1}\parens{\tfrac{4}{\alpha L}}^m \parens{\expec[]{\inorm{w_0}^m} + \parens{\tfrac{\DeltaN}{2}}^m} }\\
    &= \Delta^m a^{mt}\landau{t}{t^{m+1}}\\
    &= \Delta^m a^{mt}\landau{t}{t^{\beta+1}}, 
\end{align*}
where $q > \frac{\beta}{\eps}$ implies that $\beta > m = p(\beta-\eps)$. Note that the term $\landau{t}{t^{\beta+1}}$ is $\landau{\Delta}{1}$. Then finally we have,
\begin{equation}\label{eqstep:schemeAseries2}
    \expec[\state{}]{\inorm{y_t}^{p(\beta-\eps)}}^{\frac{1}{p}} \leq
    \Delta^{\beta-\eps}\landau{t}{t^{\frac{\beta+1}{p}}
      a^{(\beta-\eps)t}}, 
\end{equation}
which in combination with \eqref{eqstep:schemeAseries1} yields that \eqref{eqstep:schemeAmoment4} is bounded by
\begin{align*}
    \eqref{eqstep:schemeAmoment4} \leq \Delta^{\beta-\eps}\sum_{t=1}^\infty \landau{t}{t^{\beta+1}\parens{a^{\beta-\eps}\xi^{-\frac{\beta}{q}}}^t} ,
\end{align*}
where the series above converges to a constant $M>0$ when $a^{\beta-\eps}\xi^{-\frac{\beta}{q}} < 1$. This is ensured by the condition $q < \frac{\beta}{\beta-\eps}\log_a\xi$, since then we have
\begin{align*}
    a^{\beta-\eps}\xi^{-\frac{\beta}{q}} < a^{\beta-\eps}\xi^{-\beta(\frac{\beta-\eps}{\beta})\frac{1}{\log_a\xi}} = a^{\beta-\eps}\xi^{-(\beta-\eps)\log_\xi a} = a^{\beta-\eps}a^{-(\beta-\eps)} = 1.
\end{align*}
This establishes \eqref{eqstep:schemeAmomentMain} and completes the proof.
\end{proof}

\begin{proof}[Proof of Proposition \ref{prop:bin_moment_sup}]
Recall that we must show that
\begin{equation*}
    \sup_{N\geq 2}\expec[]{(\Deltapi)^{\beta - \eps}} < \infty.
\end{equation*}
We will do this again via random-time Lyapunov drift arguments. In particular, we will use most of the same Lyapunov parameters as in the previous proof. To be precise, with $r \coloneqq \rho^{\beta-\eps}$ we once again set
\begin{align*}
    &V(\state{}) = \Delta^{\beta-\eps},\\*
    &d(\state{}) = s\Delta^{\beta-\eps},\\
    &C = \Lambda\cap\cbrackets{\statep{}\in\Rn\times\Omega_\Delta : \Delta\leq D},\\
    &b = D^{\beta-\eps}\parens{\alpha^{\beta-\eps}\sup_{\state{}\in\Lambda}\expec[\state{}]{r^{\tau_\Lambda-1}}-1+s},
\end{align*}
where $s\in(0,1-\alpha^{\beta-\eps})$ is arbitrary, $D>0$ is such that
\begin{equation*}
    \statep{}\in\Lambda\textrm{ and }\Delta > D \implies \alpha^{\beta-\eps}\expec[\state{}]{r^{\tau_\Lambda-1}}-1+s \leq 0
\end{equation*}
(such a $D$ exists by Proposition \ref{prop:returntime_schemeA}), and $b$ is finite by Proposition \ref{prop:returntime_schemeA}.

Instead of $f$ proportional to $\inorm{x}^{\beta-\eps}$, we set $f$ proportional to $\Delta^{\beta-\eps}$. That is, we set:
\begin{equation*}
    f(\state{}) = cs\Delta^{\beta-\eps}, 
\end{equation*}
where
\begin{equation*}
    c = \parens{ 1 + \frac{\alpha^{\beta-\eps}}{r-1}\parens{ \sup_{\statep{}\in\Lambda}\expec[]{r^{\tau_\Lambda - 1}} - 1 } }^{-1}
\end{equation*}
is well-defined by Proposition \ref{prop:returntime_schemeA}. We will show that these Lyapunov parameters satisfy the drift condition \eqref{eq:driftcriteria_foraset} and invoke Lemma \ref{lem:drift_moment} to complete the proof.

To begin, note that the first drift inequality of \eqref{eq:driftcriteria_foraset} does not involve $f$, and with otherwise identical Lyapunov parameters we showed in the proof of Proposition \ref{prop:state_moment_sup} that this drift inequality holds. Therefore, what remains is to show only that the second inequality in \eqref{eq:driftcriteria_foraset} holds. Explicitly, we must show that
\begin{equation*}
    cs\expec[\state{}]{\sum_{t=0}^{\tau_\Lambda-1} \Delta_t^{\beta-\eps} } \leq s\Delta^{\beta-\eps},
\end{equation*}
which is equivalent to,
\begin{equation}\label{eqstep:bin_moment_mustshow}
    \expec[\state{}]{\sum_{t=0}^{\tau_\Lambda-1} \parens{\frac{\Delta_t}{\Delta}}^{\beta-\eps}} \leq c^{-1} = 1 + \frac{\alpha^{\beta-\eps}}{r-1}\parens{ \sup_{\statep{}\in\Lambda}\expec[]{r^{\tau_\Lambda - 1}} - 1 }.
\end{equation}
We establish this now. Note that $\Delta_0 = \Delta$ and for $1\leq t\leq \tau_\Lambda$ we have $\Delta_t = \alpha\Delta\rho^{t-1}$ so,
\begin{align*}
    \expec[\state{}]{\sum_{t=0}^{\tau_\Lambda-1} \parens{\frac{\Delta_t}{\Delta}}^{\beta-\eps}} &= 1 + \alpha^{\beta-\eps}\expec[\state{}]{\sum_{t=1}^{\tau_\Lambda-1} \parens{\rho^{t-1}}^{\beta-\eps}}\\
    &= 1 + \alpha^{\beta-\eps}\expec[\state{}]{\sum_{t=1}^{\tau_\Lambda-1} r^{t-1}}\\
    &= 1 + \alpha^{\beta-\eps}\cdot\frac{\expec[\state{}]{r^{\tau_\Lambda-1}} - 1}{r-1}\\
    &= 1 + \frac{\alpha^{\beta-\eps}}{r-1}\parens{\expec[\state{}]{r^{\tau_\Lambda}} - 1} \leq c^{-1}
\end{align*}
by construction of $c$. Therefore the drift condition \eqref{eq:driftcriteria_foraset} is satisfied for choice of $b$ and $f$ independent of $N$, which by Lemma \ref{lem:drift_moment} completes the proof.
\end{proof}

\bibliographystyle{IEEEtran}
%bibliography{references-jonathan,SerdarBibliography}

\end{document}